
\let\Begin\begin
\let\End\end


\documentclass[twoside,12pt]{article}
\setlength{\textheight}{24cm}
\setlength{\textwidth}{16cm}
\setlength{\oddsidemargin}{2mm}
\setlength{\evensidemargin}{2mm}
\setlength{\topmargin}{-15mm}
\parskip2mm


\usepackage[usenames,dvipsnames]{color}
\usepackage{amsmath}
\usepackage{amsthm}
\usepackage{amssymb}
\usepackage[mathcal]{euscript}

\usepackage{cite}
%

\usepackage[T1]{fontenc}
\usepackage[english]{babel}
\usepackage[babel]{csquotes}

\usepackage{latexsym}
\usepackage{graphicx}
\usepackage{mathrsfs}
\usepackage{mathrsfs}
\usepackage{hyperref}
\usepackage{pgfplots}

%
%


\definecolor{viola}{rgb}{0.3,0,0.7}
\definecolor{ciclamino}{rgb}{0.5,0,0.5}

\def\gianni #1{{\color{ciclamino}#1}}
\def\pier #1{{\color{blue}#1}}
\def\betti #1{{\color{red}#1}}
\def\juerg #1{{\color{Green}#1}}
\def\nju #1{{\color{Green}#1}}
\def\juju #1{{\color{red}#1}}

\def\pier #1{#1}
\def\nju #1{#1}
\def\betti #1{#1}
\def\juerg #1{#1}
\def\gianni #1{#1}
\def\juju #1{#1}


\bibliographystyle{plain}


%

\def\Beq{\Begin{equation}}
\def\Eeq{\End{equation}}
\def\Bsist{\Begin{eqnarray}}
\def\Esist{\End{eqnarray}}

\def\Bthm{\Begin{theorem}}
\def\Ethm{\End{theorem}}
\def\Blem{\Begin{lemma}}
\def\Elem{\End{lemma}}
\def\Bprop{\Begin{proposition}}
\def\Eprop{\End{proposition}}

\def\Brem{\Begin{remark}\rm}
\def\Erem{\End{remark}}

\def\Bdim{\Begin{proof}}
\def\Edim{\End{proof}}
\def\Bcenter{\Begin{center}}
\def\Ecenter{\End{center}}
\let\non\nonumber




\def\step #1 \par{\medskip\noindent{\bf #1.}\quad}


\def\Lip{Lip\-schitz}
\def\Holder{H\"older}
\def\Frechet{Fr\'echet}
\def\aand{\quad\hbox{and}\quad}

\def\lhs{left-hand side}
\def\rhs{right-hand side}

\def\nbh{neighborhood}
\def\bhv{behavior}

\let\hat\widehat


\def\multibold #1{\def\arg{#1}%
  \ifx\arg\pto \let\next\relax
  \else
  \def\next{\expandafter
    \def\csname #1#1#1\endcsname{{\bf #1}}%
    \multibold}%
  \fi \next}

\def\pto{.}

\def\multical #1{\def\arg{#1}%
  \ifx\arg\pto \let\next\relax
  \else
  \def\next{\expandafter
    \def\csname cal#1\endcsname{{\cal #1}}%
    \multical}%
  \fi \next}


\def\multimathop #1 {\def\arg{#1}%
  \ifx\arg\pto \let\next\relax
  \else
  \def\next{\expandafter
    \def\csname #1\endcsname{\mathop{\rm #1}\nolimits}%
    \multimathop}%
  \fi \next}

\multibold
qwertyuiopasdfghjklzxcvbnmQWERTYUIOPASDFGHJKLZXCVBNM.

\multical
QWERTYUIOPASDFGHJKLZXCVBNM.

\multimathop
diag dist div dom mean meas sign supp .

\def\Span{\mathop{\rm span}\nolimits}


\def\accorpa #1#2{\eqref{#1}--\eqref{#2}}
\def\Accorpa #1#2 #3 {\gdef #1{\eqref{#2}--\eqref{#3}}%
  \wlog{}\wlog{\string #1 -> #2 - #3}\wlog{}}


\def\separa{\noalign{\allowbreak}}

\def\somma #1#2#3{\sum_{#1=#2}^{#3}}

\def\graffe #1{\mathopen\{#1\mathclose\}}

\def\<#1>{\mathopen\langle #1\mathclose\rangle}
\def\norma #1{\mathopen \| #1\mathclose \|}

\def\[#1]{\mathopen\langle\!\langle #1\mathclose\rangle\!\rangle}

\def\iot {\int_0^t}
\def\ioT {\int_0^T}
\def\intet{\int_{Q^t}}
\def\intQt{\int_{Q_t}}
\def\intQ{\int_Q}
\def\iO{\int_\Omega}

\def\dt{\partial_t}
\def\dn{\partial_\nu}

\def\cpto{\,\cdot\,}

\def\checkmmode #1{\relax\ifmmode\hbox{#1}\else{#1}\fi}
\def\aeO{\checkmmode{a.e.\ in~$\Omega$}}
\def\aeQ{\checkmmode{a.e.\ in~$Q$}}
\def\aet{\checkmmode{a.e.\ in~$(0,T)$}}

\def\aat{\checkmmode{for a.a.~$t\in(0,T)$}}


\def\erre{{\mathbb{R}}}
\def\erren{\erre^n}

\def\enne{{\mathbb{N}}}




\def\genspazio #1#2#3#4#5{#1^{#2}(#5,#4;#3)}
\def\spazio #1#2#3{\genspazio {#1}{#2}{#3}T0}
\def\spaziot #1#2#3{\genspazio {#1}{#2}{#3}t0}

\def\L {\spazio L}
\def\H {\spazio H}

\def\Lt {\spaziot L}

\def\C #1#2{C^{#1}([0,T];#2)}


\def\Lx #1{L^{#1}(\Omega)}
\def\Hx #1{H^{#1}(\Omega)}

\def\LQ #1{L^{#1}(Q)}

\def\CQ #1{C^{#1}(\overline Q)}

\def\Luno{\Lx 1}
\def\Ldue{\Lx 2}
\def\Linfty{\Lx\infty}

\def\Huno{\Hx 1}
\def\Hdue{\Hx 2}


\def\LQ #1{L^{#1}(Q)}


\let\theta\vartheta
\let\eps\varepsilon
\let\phi\varphi

\let\TeXchi\chi                         
\newbox\chibox
\setbox0 \hbox{\mathsurround0pt $\TeXchi$}
\setbox\chibox \hbox{\raise\dp0 \box 0 }
\def\chi{\copy\chibox}


\def\QED{\hfill $\square$}


\def\CO{C_\Omega}
\def\Vp{V^*}
\def\Vn{V_n}

\def\phieps{\phi^\eps}
\def\mueps{\mu^\eps}

\def\veps{v_\eps}
\def\phiepsk{\phieps_k}
\def\Psieps{\Psi_\eps}
\def\geps{g_\eps}

\def\Beta{\hat\beta}
\def\betaeps{\beta_\eps}
\def\Betaeps{\hat\beta_\eps}
\def\betaz{\beta^\circ}
\def\Pi{\hat\pi}

\def\phiz{\phi_0}

\def\phin{\phi_n}
\def\mun{\mu_n}
\def\xin{\xi_n}
\def\phinj{\phi_n^j}
\def\munj{\mu_n^j}

\def\barphi{\overline \phi}

\def\baru{\overline u}
\def\barv{\overline v}
\def\vstar{v_*}
\def\barvstar{\barv_*}

\def\barphin{\overline\phin}
\def\barmun{\overline\mun}
\def\barphiz{\barphi_0}

\def\zeps{z_\eps}
\def\zk{\zeta_k}

\def\pn{p_n}
\def\qn{q_n}
\def\barpn{\overline\pn}

\def\phiQ{\phi_Q}
\def\phiO{\phi_\Omega}
\def\muQ{\mu_Q}

\def\S{{\cal S}}
\def\Uad{\calU_{ad}}
\def\UR{\calU_R}
\def\ustar{u^*}
\def\phistar{\phi^*}
\def\mustar{\mu^*}

\def\normaV #1{\norma{#1}_V}
\def\normaW #1{\norma{#1}_W}
\def\normaVp #1{\norma{#1}_*}

\Begin{document}

%
\title{\gianni{Well-posedness and optimal control\\ for a Cahn--Hilliard--Oono system\\
with control in the mass term}}
\author{}
\date{}
\maketitle
\Bcenter
\vskip-1.5cm
{\large\sc Pierluigi Colli$^{(1)}$}\\
{\normalsize e-mail: {\tt pierluigi.colli@unipv.it}}\\[.25cm]
{\large\sc Gianni Gilardi$^{(1)}$}\\
{\normalsize e-mail: {\tt gianni.gilardi@unipv.it}}\\[.25cm]
{\large\sc Elisabetta Rocca$^{(1)}$}\\
{\normalsize e-mail: {\tt elisabetta.rocca@unipv.it}}\\[.25cm]
{\large\sc J\"urgen Sprekels$^{(2)}$}\\
{\normalsize e-mail: {\tt sprekels@wias-berlin.de}}\\[.45cm]
$^{(1)}$
{\small Dipartimento di Matematica ``F. Casorati'', Universit\`a di Pavia}\\
{\small and Research Associate at the IMATI -- C.N.R. Pavia}\\
{\small via Ferrata 5, 27100 Pavia, Italy}\\[.2cm]
$^{(2)}$
{\small Department of Mathematics}\\
{\small Humboldt-Universit\"at zu Berlin}\\
{\small Unter den Linden 6, 10099 Berlin, Germany}\\[2mm]
{\small and}\\[2mm]
{\small Weierstrass Institute for Applied Analysis and Stochastics}\\
{\small Mohrenstrasse 39, 10117 Berlin, Germany}
\Ecenter
\begin{center}
\emph{\gianni{Dedicated to our dear friend Maurizio Grasselli\\
on the occasion of his 60th birthday}}
\end{center}
\Begin{abstract}
\betti{The paper treats the problem of optimal distributed control of a \pier{Cahn--Hilliard--Oono} 
system in $\erre^d$, $1\leq d\leq 3$, \pier{with the control located in 
the mass term and admitting} general potentials that include both the case of \pier{a} 
regular potential and the case of \pier{some} singular potential. 
The first part of the paper \pier{is concerned with} the dependence of the phase variable 
\pier{on} the control variable. For this purpose, suitable regularity and continuous 
dependence results are shown. In particular, in \pier{the} case of a logarithmic 
potential, we need to prove an ad hoc strict separation property, and for this reason we 
\pier{have to restrict ourselves} to the case $d=2$. \pier{In the rest of the work, we study  the necessary first-order optimality conditions, which} are proved under suitable compatibility conditions on the initial datum of the phase \pier{variable and the time derivative of the} control, at least in case of potentials having \pier{unbounded} domain.}

\End{abstract}

\pagestyle{myheadings}
\newcommand\testopari{\sc Colli \ --- \ Gilardi \ --- \ Rocca \ --- \ Sprekels}
\newcommand\testodispari{\sc \gianni{The Cahn--Hilliard--Oono system
with control in the mass term}}
\markboth{\testopari}{\testodispari}


\section{Introduction}
\label{Intro}
\setcounter{equation}{0}

\betti{In this paper, we study the following PDE system, referred \juju{to} in the literature as 
\pier{the Cahn--Hilliard--Oono} (CHO) system (cf.~\cite{oonopuri87, oonopuri88I, oonopuri88II}), which is of great interest in the study of pattern formations in 
phase-separating materials}\pier{:} 
\Bsist
  & \dt\phi + \epsilon(\phi-\hat c\, ) - \Delta\mu = 0
  & \quad \hbox{in $Q:=\Omega\times(0,T)$},
  \label{primaCHO}
  \\
  & \mu = - \Delta {\phi + f'(\phi)} 
  & \quad \hbox{in $Q$},
  \label{secondaCHO}
  \\
  & \dn\mu = \dn\phi = 0 
  & \quad \hbox{on $\partial\Omega\times(0,T)$},
		\label{Ibc} 
  \\
  & \phi(0) = \phiz
  & \quad \hbox{in $\Omega$}.
  \label{Icauchy}
\Esist
\juju{Here,} $\Omega$ is the domain \betti{in $\erre^d$, $1\leq d\leq 3$,}  where the evolution takes place, and $T$ is some final time. \pier{In \eqref{Ibc},} $\dn$ denotes the outward normal derivative on the boundary \pier{$\partial \Omega$}. In the above equations, 
the unknowns are~$\phi$, the order parameter \betti{representing the relative monomer concentration difference}, and~$\mu$, the chemical potential, 
while $\hat c$~is a prescribed mass average, $\epsilon$~is a given phenomenological constant, 
and $f'$~is the derivative of a double-well potential~$f$.
Finally, $\phiz$ is a given initial datum. 

Typical and important examples for $f$ 
are the so-called {\em classical regular potential} and the {\em logarithmic double-well potential\/}{,
which are the semiconvex functions given~by}
\Bsist
  && f_{reg}(r) := \frac 14 \, (r^2-1)^2 \,,
  \quad r \in \erre ,
  \label{regpot}
  \\
  && f_{log}(r) := (1+r)\ln (1+r)+(1-r)\ln (1-r) - c_1 r^2 \,,
  \quad r \in (-1,1),
  \label{logpot}
\Esist
where $c_1>1$ so that $f_{log}$ is nonconvex.
Another example is {the {\em double obstacle potential\/}, where, with} $c_2>0$,
\Beq
  f_{2obs}(r) := - c_2 r^2 
  \quad \hbox{if $|r|\leq1$}
  \aand
  f_{2obs}(r) := +\infty
  \quad \hbox{if $|r|>1$}.
  \label{obspot}
\Eeq
In cases like \eqref{obspot}, one has to split $f$ into a nondifferentiable convex part 
(the~indicator function of $[-1,1]$ in the present example) and a smooth perturbation.
Accordingly, one has to replace the derivative of the convex part
by the subdifferential and read \eqref{secondaCHO} as a differential inclusion. \betti{The coefficients $c_1$ and $c_2$ are related to the \juju{critical absolute} temperature at which the phase separation takes place.}

\betti{System \accorpa{primaCHO}{secondaCHO} can be seen as a Cahn\pier{--}Hilliard equation with reaction and 
\juju{turns} out to be useful in several applications 
such as biological models~\cite{KS}, inpainting algorithms~\cite{BEG}, and 
polymers~\cite{BO}. Indeed, taking $\epsilon=0$, the CHO equation reduces to the classical Cahn\pier{--}Hilliard (CH) equation (see \cite{BE,C,CH}). The main feature of \eqref{primaCHO} is that, contrary to the classical CH equation, it does not imply the \juju{conservation of total mass}. Indeed, integrating it over $\Omega$ and using the boundary \pier{and initial conditions \eqref{Ibc}--\eqref{Icauchy}, one finds \juju{that}}
$$
\barphi(t)=\hat c +{\rm e}^{-\epsilon t}(\barphiz -\hat c) \quad\hbox{for any }\,  \pier{t\geq 0},
$$
where \pier{$\barphi (t) :=\frac{1}{|\Omega|}\int_{\Omega} \varphi(t)\, dx$ denotes the mean value over $\Omega$ at time $t$.}
We notice that we have mass conservation only in \juju{the}
case $\barphiz=\hat c$. In all other cases, \juju{the time function $\barphi(t)$ converges exponentially fast to $\hat c$
as $t $ goes to  $+ \infty$}. 
Hence, the \pier{reaction term present in \eqref{primaCHO} accounts for long-range interactions. Of course, more general terms, in particular nonlinearities, could be considered} (cf., 
for example, \cite{MR}, where the nonlocal CH equation with nonlinear reaction term was analyzed).
Moreover, let us notice that the interest of \pier{including} reaction terms in CH-type \pier{systems} is particularly increasing \pier{due to the development of diffuse interface models of tumor growth coupling CH with nutrient diffusion and other equations. Indeed, in this class of models the parameter $\varphi$ should be intended as the concentration of the tumor phase, and so, in this case, it is of particular interest not to have \juju{mass} conservation, but to observe a possible growth or decrease of the \juju{tumor mass}, due to \juju{the} proliferation or death of cells. In this respect, we may quote \cite{CRW, CGH, CGRS, CGRS1, CGRS2, CGS-2019a, FGR, GLR, GARL_1, GARL_3, MRS, S_b, S_a} and \juju{the}
references therein for results on well-posedness, long-time \juju{behavior} of solutions, \pier{asymptotic analyses}, and optimal control, related to tumor growth models.} More recently, always in the framework of tumor growth dynamics, \pier{in~\cite{GLRS}} a generalized CHO equation has been coupled with an Hele--Shaw equation, and \juju{global} 
well-posedness and regularity results have been proved in the 2D case.}

\betti{The \pier{solvability} and the existence of global and exponential attractors for the CHO equation with regular potential \eqref{regpot} have been studied in \cite{M}, while in \cite{GGM} the authors \juju{investigated}
 the case of the singular logarithmic type potential \eqref{logpot}\pier{, by} establishing some regularization properties of the unique solution in finite time. Both in 2D and 3D they \juju{proved} the existence of a global attractor; moreover, in the 2D case, taking advantage of their proof of a strict separation property, \pier{the authors of \cite{GGM} also \juju{showed}} the existence of \juju{an} exponential attractor and the convergence to a single equilibrium. 
Finally, the CHO equation has been coupled with Navier--Stokes systems both in case of regular \juju{potentials}
~\cite{BGM} and of singular \juju{potentials}~\cite{MT}.}

In the present work, we study a distributed control problem for the \pier{system \accorpa{primaCHO}{Icauchy}, where the term $\epsilon(\phi-\hat c\, )$ in 
\eqref{primaCHO} is replaced by $\phi- u$, thus taking $\epsilon = 1 $ and substituting $\hat c$ \juerg{by}
 a control $u$ which is allowed to be a function of space and time.}
\betti{Having in mind the application to tumor growth models, the control $u$  could represent a source of therapy, like, for example, a concentration of \juju{cytotoxic} drugs introduced in the cells in order to minimize the volume of the tumor (\juju{which is} represented by the integral of $\varphi$ over the domain $\Omega$) or to reach some target tumor distribution at the final time $T$ of the therapy cycle.} \pier{About optimal control problems for CH systems, let us mention some related work. A very general approach for distributed control problems for possibly fractional equations of CH-type is carried out in the papers \cite{CGS-2019c, CGS-2020, CGS-2021}, with \juju{an} extension of the analysis to double obstacle potentials like $ f_{2obs}$ in \eqref{obspot} via deep quench approximation. The coupling of CH equations in the bulk with dynamic boundary conditions has been investigated in \cite{CGS-2014, CGS-2015, CS-2021},  and the presence of a convective term with the velocity vector taken as control has been dealt with in \cite{CGS-2018a, CGS-2018b, CGS-2019b, GS-2019} (see also \juju{the} references in the quoted contributions).}

\pier{We now go back to the goal of this paper and observe that, \juju{since we are not interested in the longtime \bhv\ of the solution, the constants $\hat c$ and $\epsilon$ here actually} do not play any significant role, whence} we can assume that $\epsilon=1$ and absorb $\hat c$ in~$\,u$.
Therefore, the state system under investigation is the following:
\Bsist
  & \dt\phi + \phi - \Delta\mu = u
  & \quad \hbox{in \,$Q$,}
  \label{Iprima}
  \\
  & \mu = - \Delta{\phi + f'(\phi)} 
  & \quad \hbox{in \,$Q$,}
  \label{Iseconda}
\Esist
with the same boundary and initial conditions as before.
For our purpose, it is crucial to study the dependence of the solution on the control variable~$u$.
Hence, a major part of the paper is devoted to the well-posedness of the system
and the continuous dependence just mentioned.
It must be pointed out that, at least in the case of potentials having a bounded domain,
the control $\,u\,$ and the initial datum $\,\phiz\,$ must satisfy 
proper necessary compatibility conditions {in order to guarantee} the existence of a solution.
For this reason, even well-posedness cannot be deduced from \cite{GGM} and the existing literature.
Then, we study the control problem. 
Precisely, we want to minimize the tracking-type cost functional
\Beq
  \calJ({(\phi,\mu)},u)
  := \frac{\alpha_1}2 \intQ |\phi - \phiQ|^2
  + \frac{\alpha_2}2 \iO |\phi(T) - \phiO|^2
  + \frac{\alpha_3}2 \intQ |\mu - \muQ|^2
  + \frac{\alpha_4}2 \intQ |u|^2
  \label{Icost}
\Eeq
over the set of admissible controls
\Beq
  \Uad := \graffe{u \in \H1H\cap\LQ\infty: \ \norma u_{\LQ\infty} \leq M, \ \norma{\dt u}_{\LQ2} \leq M'},
  \label{IUad}
\Eeq
subject to the system given by \accorpa{Iprima}{Iseconda} and \accorpa{Ibc}{Icauchy}.
In \eqref{Icost}, $\phiQ$ and $\muQ$ are given functions on~$Q$, 
$\phiO$~is a given function on $\Omega$, and $\alpha_i$, $i=1,\dots,4$, are nonnegative constants (not all zero).
In \eqref{IUad}, $M>0$ and $M'>0$ are prescribed constants as well.
We \pier{show} the existence of an optimal control.
Then the main point is {to prove} the \Frechet\ differentiability 
of the control-to-state operator between suitable functional spaces,
which allows us to establish first-order necessary optimality conditions 
in terms of the solution to the linearized system; the latter is then eliminated 
by means of the solution to a proper adjoint problem.

The paper is organized as follows. 
In the next section, we list our assumptions and notations
and state our results.
The proofs of 
those regarding the well-posedness of the problem, 
the continuous dependence of its solution on the control variable, and the regularity,
are given in Sections~\ref{CONTDEP}--\ref{REGULARITY}.
The argument for the existence and the regularity of the solution
is based on the study of proper approximating problems 
performed in the first part of Section~\ref{EXISTENCE}.
Finally, the last \pier{Section~\ref{CONTROL} is devoted to the analysis of} 
the control problem and the corresponding first-order necessary optimality conditions.


\section{Statement of the problem and results}
\label{STATEMENT}
\setcounter{equation}{0}

In this section, we state precise assumptions and notations and present our results.
First of all, the subset $\Omega\subset\erre^d$ with $1\leq d\leq3$ is assumed to be bounded, connected and smooth.
As in the Introduction,   
$\dn$ stands for the normal derivative on $\Gamma:=\partial\Omega$.
Moreover, we set for brevity
\Beq
  Q_t := \Omega \times (0,t)
  \aand
  Q^t := \Omega \times (t,T)
  \quad \hbox{for $0<t<T$},
  \aand
  Q := \Omega \times (0,T).
  \label{defQt}
\Eeq
If $X$ is a Banach space, $\norma\cpto_X$ denotes both its norm and the norm of~$X^d$.
The only exception from this convention on the norms is given
by the spaces $L^p$ ($1\leq p\leq\infty$) constructed on $(0,T)$,
$\Omega$ and~$Q$, 
whose norms are often denoted by~$\norma\cpto_p$,
and by the space $H$ defined below,  
whose norm is simply denoted by~$\norma\cpto$.
We~put
\Bsist
  && H := \Ldue \,, \quad  
  V := \Huno 
  \aand
  W := \graffe{v \in \Hdue: \ \dn v = 0} .
  \label{defspazi}
\Esist 
Moreover, $\Vp$~is the dual space of $V$ and $\<\cpto,\cpto>$ is the dual pairing between $\Vp$ and~$V$.
In the following, we work in the framework of the Hilbert triplet
$(V,H,\Vp)$.
Thus, by also using the symbol $(\cpto,\cpto)$ for the standard inner product of $H$,
we have
$\<g,v>=(g,v)$
for every $g\in H$ and $v\in V$. {We also use the symbol  
$(\cpto,\cpto)$ for the standard inner product in any of the product spaces $H^N$ for $N\in\enne$.}

Next, we introduce the generalized mean value.
We write $|\Omega|$ for the measure of~$\Omega$ and~set
\Beq
  \barvstar := |\Omega|^{-1} \< \vstar,1 >
  \quad \hbox{for $\vstar\in\Vp$}
  \label{defmean}
\Eeq
where $1$ denotes the constant function $x\mapsto1$, $x\in\Omega$.
More generally, to simplify the notation, we use the same symbol 
for the real number $a$ and the associated constant functions on $\Omega$ and~$Q$.
It is clear that $\barvstar$ is the usual mean value of $\vstar$ if $\vstar\in H$.

Now, we list our assumptions on the structure of our system and the data. 
For the potential $\,f$, we assume that
\Bsist
  && f: \erre \to (-\infty,+\infty]
  \quad \hbox{can be split as} \quad
  f = \Beta + \Pi,
  \quad \hbox{where}
  \label{hpf}
  \\
  && \Beta: \erre \to [0,+\infty]
  \quad \hbox{is convex, proper, and l.s.c.\ with $\Beta(0)=0$,}
  \label{hpBeta}
  \\
  && \Pi: \erre \to \erre 
  \quad \hbox{is of class $C^1$, and its derivative is \Lip\ continuous}.
  \qquad
  \label{hpPi}
\Esist
\Accorpa\HPstruttura hpf hpPi
We set for convenience
\Beq
  \beta := \partial\Beta 
  \aand
  \pi := \Pi \, {}',
  \label{defbetapi}  
\Eeq
and notice that $\beta$ is a maximal monotone graph in $\erre\times\erre$ 
satisfying $0\in\beta(0)$.
We use the symbols $D(\beta)$ and $\betaz(r)$ 
for the domain of $\beta$ and 
the element of $\beta(r)$ (with $r\in D(\beta)$) having minimum modulus.
We extend the notations $\Beta$, $\beta$, $D(\beta)$ and $\betaz$
to the functionals and the operators induced on $L^2$ spaces.

Even though the control variable $u$ is fixed when dealing with well-posedness,
we prepare the possibility of letting $u$ vary by introducing an upper bound $M$ in the data.
We assume that
\Bsist
  && u \in \LQ\infty , \quad
  M \in [0,+\infty)
  \aand
  \norma u_\infty \leq M,
  \label{hpu}
  \\
  && \phiz \in W
  \quad \hbox{with $\,\inf\phiz\,$ and $\,\sup\phiz\,$ belonging to the interior of $D(\beta)$},
  \qquad
  \label{hpphiz}
  \\
  && \barphiz \pm M
  \quad \hbox{belong to the interior of $D(\beta)$}.
  \label{compatibility}
\Esist
\Accorpa\HPdati hpu compatibility
Assumption \eqref{hpphiz} is rather strong{; its first application} 
and a comment on the last assumption \eqref{compatibility} are given 
in the forthcoming Remarks~\ref{Phinz} and~\ref{Compatibility}, respectively.

Let us come to the definition of {our notion of} solution, which we state in a weak form.
Namely, a solution is a triplet $(\phi,\mu,\xi)$ satisfying the regularity requirements
\Bsist
  && \phi \in \H1\Vp \cap \L\infty V,
  \label{regphi}
  \\
  && \mu \in \L2V,
  \label{regmu}
  \\
  && \xi \in \L2H,
  \aand
  \xi \in \beta(\phi) \quad \aeQ ,
  \label{regxi}
\Esist
\Accorpa\Regsoluz regphi regxi
and the following variational equations and initial condition:
\Bsist
  && \< \dt\phi(t) , v > 
  + \iO \phi(t) v 
  + \iO \nabla\mu(t) \cdot \nabla v
  = \iO u(t) v
  \non
  \\
  && \quad \hbox{\aat\ and every $v\in V$},
  \label{prima}
  \\[1mm]
  && \iO \mu(t) v 
  = \iO \nabla\phi(t) \cdot \nabla v
  + \iO \xi(t) v
  + \iO \pi(\phi(t)) v
  \non
  \\
  && \quad \hbox{\aat\ and every $v\in V$},
  \label{seconda}
  \\[2mm]
  && \phi(0) = \phiz \quad \hbox{\betti{a.e.~in }}\betti{\Omega} \,.
  \label{cauchy}
\Esist
\Accorpa\Pbl prima cauchy

\Brem
\label{IntPbl}
We observe that the above variational equations are equivalent 
to their time-integrated versions with time dependent test functions,~i.e.,
\Bsist
  && \ioT \< \dt\phi(t) , v(t) > \, dt
  + \intQ \phi \, v
  + \intQ \nabla\phi \cdot \nabla v
  = \intQ u v,
  \label{intprima}
  \\
  && \intQ \mu v 
  = \intQ \nabla\phi \cdot \nabla v
  + \intQ \xi v
  + \iO \pi(\phi) v,
  \label{intseconda}
\Esist
both for every $v\in\L2V$.
We also point out that \eqref{seconda} can be written as a boundary value problem.
Namely, since $\phi\in\L2V$ and $\mu-\xi-\pi(\phi)\in\L2H$,
it is equivalent~to
\Beq
  \phi \in \L2W,
  \aand
  \mu = -\Delta\phi + \xi + \pi(\phi)
  \quad \aeQ ,
  \label{secondabvp}
\Eeq
by elliptic regularity.
On the contrary, the analogue for \eqref{prima} would be true 
only if $\dt\phi$ {were} more regular.
\Erem

\Brem
\label{Compatibility}
Let us comment on assumption \eqref{compatibility}.
By just taking $v=1/|\Omega|$ in \eqref{prima} and accounting for \eqref{cauchy}, we have that
\Beq
  \frac d{dt} \, \barphi(t) + \barphi(t) = \baru(t)
  \quad \aat,
  \aand
  \barphi(0) = \barphiz,
  \label{ode}
\Eeq
which yield
\Beq
  \barphi(t) 
  = \barphiz + \iot e^{-(t-s)} \baru(s) \, ds
  \quad \hbox{for every $t\in[0,T]$}.
  \label{soluzode}
\Eeq
Since this implies that
\Beq
  |\barphi(t) - \barphiz|
  \leq \norma\baru_\infty
  \leq \norma u_\infty
  \leq M,
  \non
\Eeq
assumption \eqref{compatibility} ensures that $\barphi(t)$ belongs to $D(\beta)$ for every $t\in[0,T]$.
We observe that the \rhs\ of \eqref{soluzode} just depends on $\phiz$ and $u$,
and a necessary condition for the existence of a solution is that it remains in $D(\beta)$ for all times.
\gianni{In particular, if $u=2$ and $\phiz=0$, 
it is given by $2(1-e^{-t})$ and becomes larger than 1 if $t>\ln2$,
so that no solution can exist on an interval $[0,T]$ with $T>\ln2$ 
if $f$ is the logarithmic potential~\eqref{logpot}.} 

\Erem

We have the following well-posedness result.
\Bthm
\label{Wellposedness}
Suppose that the conditions \HPstruttura\ on the structure of the system and \HPdati\ on the data are fulfilled.
Then there exists \gianni{at least one triplet $(\phi,\mu,\xi)$,
with the regularity conditions \Regsoluz, that solves problem \Pbl\
and satisfies the estimate}
\Beq
  \norma\phi_{\H1\Vp\cap\L\infty V\cap\L2W}
  + \norma\mu_{\L2V}
  + \norma\xi_{\L2H}
  \leq K_1
  \label{stimaK1}
\Eeq
with a constant $K_1>0$ that depends only on the structure of the system, 
$\Omega$, $T$, the initial datum, and~$M$. 

Furthermore, if {$u_i\in\LQ\infty$, $i=1,2$, are given} 
and $(\phi_i,\mu_i,\xi_i)$ are \gianni{any} corresponding solutions, then
the estimate
\Beq
  \norma{\phi_1-\phi_2}_{\C0\Vp\cap\L2V}
  \leq C_1 \, \bigl( \norma{u_1-u_2}_{\L1\Vp} + \norma{\baru_1-\baru_2}_{L^1(0,T)}^{1/2} \bigr)
  \label{contdep1}
\Eeq
holds true with a constant $C_1>0$ that depends only on the structure of the system, 
$\Omega$, $T$, and \gianni{an upper bound for $\norma{\xi_i}_{\pier{L^1(Q)}}$, $i=1,2$}. 
{In particular, the solution component $\,\phi\,$
is uniquely determined. If, in addition, $\,\beta\,$ is single-valued, then also $\,\mu\,$ and $\,\xi\,$
are uniquely determined.} 
\Ethm

An improvement of the regularity of the solution can be {achieved} under stronger assumptions on the data.
Namely, we also assume:
\Bsist
  && \vbox{\hsize0.6\hsize\noindent \gianni{the interior of $D(\beta)$ contains $0$ \hfill\break\indent
  and the restriction of $\beta$ to it is a $C^1$ function,}}
  \label{hpbetareg}
  \\
  && u \in \H1H , \quad
  M' \in [0,+\infty)
  \aand
  \norma{\dt u}_{\L2H} \leq M' \,,
  \label{hpureg}
  \\
  && \phiz \in \Hx3 \,.
  \label{hpphizreg}
\Esist
\Accorpa\HPdatireg hpbetareg hpphizreg
\gianni{Notice that \eqref{hpbetareg} does not imply that $\beta$ is single-valued
\pier{(cf.~\eqref{obspot})}
so that uniqueness for the solution is not guaranteed.
However, we can prove the existence of a smoother solution.
Of course, if in addition $\beta$ is single-valued, then
this regularity is ensured for the unique solution.}

\Bthm
\label{Regularity}
In addition to the assumptions of Theorem~\ref{Wellposedness},
assume that \HPdatireg\ are fulfilled.
Then, \gianni{there exists at least \nju{one} solution $(\phi,\mu,\xi)$ to \Pbl\ 
that also satisfies both the regularity requirements
\Bsist
  && \phi \in \H1V \cap \L\infty W \cap \CQ0, 
  \non
  \\
  && \quad \mu \in \L\infty V
  \aand
  \xi \in \L\infty H\,, 
  \label{regularity}
\Esist
and the estimate
\Bsist
  \norma\phi_{\H1V\cap\L\infty W\cap C^0(\overline Q)}
  + \norma\mu_{\L\infty V}
  + \norma\xi_{\L\infty H}
  \leq K_2\,,
  \label{regestimate}
\Esist
with a constant $K_2$ that depends only on the structure of the system, 
$\Omega$, $T$, the initial datum, and the upper bounds $M$ and $M'$ 
appearing in \eqref{hpu} and~\eqref{hpureg}}.
\Ethm

In order to approach the control problem with a possibly singular potential, 
it is important to ensure that the solution component $\phi$ takes its values 
far away from the singularities of~$\beta$.
If $\beta$ is smooth in the interior of its domain, \nju{then} it is sufficient that all
values of $\phi$ belong to a compact interval
contained in the interior of~$D(\beta)$.
In the case of an everywhere defined smooth potential like \eqref{regpot},
the boundedness of $\phi$ given by Theorem~\ref{Regularity} is sufficient.
On the contrary, if $D(\beta)$ is a bounded open interval 
as in the case of the logarithmic potential~\eqref{logpot}, then
the requested property of $\phi$  
is equivalent to the boundedness of the solution component $\xi$.
This is related to the boundedness of~$\mu$.

\Bprop
\label{Separation}
Under the assumptions of Theorem~\ref{Wellposedness}, 
let $(\phi,\mu,\xi)$ be \gianni{a solution} to problem \Pbl.
Then boundedness of $\phi$ and $\mu$ implies boundedness of $\xi$ and the estimate
\Beq
  \norma\xi_\infty \leq \norma{\gianni{\phi+\mu-\pi(\phi)}}_\infty \,.
  \label{xibdd}
\Eeq
In particular, if $D(\beta)$ is a bounded open interval,
then all of the values of $\phi$ {are attained in} a compact interval contained in $D(\beta)$
that depends only on the shape of~$\beta$ and an upper bound for the norm of $\xi$ in $\LQ\infty$.
\Eprop

Boundedness of $\mu$ is a consequence of \eqref{regularity} if $d=1$,
since $V\subset\Linfty$ in this case.
If $d\in\{2,3\}$, the problem is open, in general.
A~natural way to ensure that $\mu$ is bounded is to prove that $\dt\phi$ belongs to $\L\infty H$.
Indeed, then \eqref{prima} {and elliptic regularity} would yield that $\mu$ belongs to $\L\infty W$, 
and thus to $\LQ\infty$ since $W\subset\Linfty$.
We cannot deal with the general case.
However, we can prove this result in the two-dimensional case if $f$ is the logarithmic potential.

\Bprop
\label{Separation2d}
Assume that $d=2$, that $f$ is the logarithmic potential \eqref{logpot},
and that the data satisfy \HPdati, \HPdatireg, as~well~as
\Beq
  \phiz \in \Hx4
  \aand
  \dn \Delta\phiz = 0 {\quad\mbox{on $\,\Gamma$}}\,.
  \label{newhpphiz}
\Eeq
Then the solution $(\phi,\mu,\xi)$ to problem \Pbl\ also satisfies
\Beq
  \dt\phi \in \L\infty H \cap \L2W \,, \quad
  \mu \in \LQ\infty \,, \quad
  \xi \in \LQ\infty\,,
  \label{reg2d}
\Eeq
as well as
\Beq
  \norma{\dt\phi}_{\L\infty H\cap\L2W}
  + \norma\mu_\infty
  + \norma\xi_\infty
  \leq K_3\,,
  \label{separation2d}
\Eeq
with a constant $K_3$ that depends only on 
$\Omega$, $T$, the initial datum, and the upper bounds $M$ and $M'$ 
appearing in \eqref{hpu} and~\eqref{hpureg}.
\Eprop

\Brem
\label{Rem2d}
We point out that under the assumptions of Proposition~\ref{Separation2d} 
the second part of Proposition~\ref{Separation}
provides a compact interval $[a,b]\subset(-1,1)$, 
which contains all of the values of $\phi$
and depends only on $\Omega$, $T$ and the constants $M$ and~$M'$.
\Erem

Whenever the potential is smooth in its domain 
and all of the values of $\phi$ belong to a compact interval $[a,b]$ contained in the interior of~$D(\beta)$,
then {the potential $f$ can be replaced a posteriori} by a potential having a \Lip\ continuous derivative.
In {this} situation, we can prove a second continuous dependence result.

\Bthm
\label{Lipcontdep}
Assume that the potential $f$ has a \Lip\ continuous derivative.
If $u_i{\in\LQ\infty}$, $i=1,2$, are two choices of $u$ 
and $(\phi_i,\mu_i,\xi_i)$ are the corresponding solutions, then
the estimate
\Beq
  \norma{\phi_1-\phi_2}_{\C0H\cap\L2W}
  + \norma{\mu_1-\mu_2}_{\L2H}
  \leq C_2 \, \norma{u_1-u_2}_{\L2H}
  \label{contdep3}
\Eeq
holds true with a constant $C_2$ that depends only on 
$\Omega$, $T$, and the \Lip\ constant of~$f'$.
\Ethm

The above results prepare the study of the control problem,
which consists in minimizing the cost functional \eqref{Icost}
over the set \eqref{IUad} of the admissible controls subject to the state system,
as said in the Introduction.
We recall the definitions and \pier{make precise assumptions}.
The cost functional is given~by
\Beq
  \calJ({(\phi,\mu)},u)
  := \frac{\alpha_1}2 \intQ |\phi - \phiQ|^2
  + \frac{\alpha_2}2 \iO |\phi(T) - \phiO|^2
  + \frac{\alpha_3}2 \intQ |\mu - \muQ|^2
  + \frac{\alpha_4}2 \intQ |u|^2\,,
  \label{cost}
\Eeq
where
\begin{align}
  &\phiQ ,\, \muQ \in \LQ2 \,, \quad
  \phiO \in \Ldue\,,
  \non\\
  &\quad{} 
  \pier{\hbox{and}\quad
  \alpha_i \in [0,+\infty)\quad \hbox{are \betti{ not all zero},} 
  \quad \hbox{for $i=1,\dots,4$} \,.}
  \label{hpcost}
\end{align}
The set of the admissible controls is given~by
\Beq
  \Uad := \graffe{u \in \H1H\cap\LQ\infty: \ \norma u_{\LQ\infty} \leq M, \ \norma{\dt u}_{\LQ2} \leq M'}\,,
  \label{defUad}
\Eeq
where
\Beq
  M ,\, M' \in [0,+\infty)
  \ \ \hbox{and \eqref{compatibility} is satisfied}.
  \label{hpUad}
\Eeq
\Accorpa\HPcontrol cost hpUad
The control problem is the following:
\Bsist
  && \vbox{\hsize .8\hsize \sl \rightskip 0pt plus 1fil \noindent
  Minimize the cost functional \eqref{cost} over the set \eqref{defUad} of admissible controls subject to the state system \Pbl.}
  \label{controlPbl}
\Esist

{At this point, it is worth noticing that the cost functional in the form \eqref{cost} is well defined
only if both solution components $\,\phi,\mu\,$ are uniquely determined. Under the assumptions of Theorem 2.3, 
this can only be guaranteed if $\,\beta\,$ is single-valued; otherwise one has to simplify the cost functional 
by postulating that $\alpha_3=0$. 
In the situation of Theorem~\gianni{\ref{Lipcontdep}}, however, this problem does not arise.   
We have the following existence result.}

\Bthm
\label{Optimum}
{Suppose that the conditions \HPstruttura, \eqref{hpphiz}, and \HPcontrol\ are fulfilled, 
and assume that
either $\alpha_3=0$ or $\beta$ is single-valued.}
Then the control problem \eqref{controlPbl} has at least {one} solution~$\ustar$.
\Ethm

The next effort is to find first-order necessary conditions for optimality.
To this end, we have to enlarge $\Uad$ in a proper topology and introduce
the control-to-state operator~$\calS$.
We fix $R>0$ small enough in order~that
\Beq
  \barphiz \pm (M+R),
  \quad \hbox{belong to the interior of $D(\beta)$},
  \label{hpR}
\Eeq
and we set
\Beq
  \UR := \graffe{u \in \H1H\cap\LQ\infty: \ \norma u_{\LQ\infty} < M+R, \ \norma{\dt u}_{\LQ2} < M'+R}.
  \label{defUR}
\Eeq
\gianni{This is an open \nbh\ of $\Uad$ in the topology 
of the first of the Banach spaces we introduce for later use:}
\Bsist
  && \calX := \H1H \cap \LQ\infty 
  \non
  \\
  && \aand
  \calY := \bigl( \C0H \cap \L2W \bigr) \times \L2H.
  \label{defXY}
\Esist
In order to continue, we need the functional $f$ to be more regular
on the range of the first component $\phi$ of the state corresponding to any element $u\in\UR$.
We first observe that {thanks to \eqref{hpR} the results stated above hold true}
with $M$ and $M'$ replaced by $M+R$ and $M'+R$, respectively; {in addition, they}
provide estimates that are uniform with respect to $u\in\UR$.
Next, we notice that in some situations it is ensured {that there exists} a compact interval $[a,b]$ 
contained in the interior of the domain of the double well potential $f$
such that, for every $u\in\UR$, the first component $\phi$ of the corresponding state
{attains its} values in~$[a,b]$.
Namely, as we have seen, {this is the case either if, under the assumptions of Theorem~\gianni{\ref{Regularity}}, 
$f$~is an everywhere defined smooth potential in any dimension {$d\in\{1,2,3\}$}
or \gianni{if the assumptions of Proposition~\ref{Separation2d} are fulfilled
(in~particular, $d=2$ and $f$ is the logarithmic potential~\eqref{logpot})}.}

Therefore, we take the following starting point:
there exist a compact interval $[a,b]$ and a constant $K>0$ such~that:
\begin{align}
  & \hbox{$f$ is a function of class $C^3$ in a \nbh\ of $[a,b]$}\pier{;}
  \label{hpfC3}
  \\[3mm]
  & {\mbox{\pier{it} holds }\,a \leq \phi \leq b
  \quad \mbox{\,a.e. in }\,Q\,\mbox{ and }\,
  \max_{0\leq i\leq 3}\norma{f^{(i)}(\phi)}_\infty \leq K,
  \hbox{ whenever $(\phi,\mu,\xi)$, }}\non\\
	&{\mbox{with }\,\xi=\beta(\phi), \mbox{ is a solution to \Pbl\ for some $u\in\UR$}}.
  \label{globalbdd}
\end{align}
\Accorpa\Start hpfC3 globalbdd
{Notice that under the conditions \Start\ also the solution component $\,\mu\,$ 
is uniquely determined so that the control-to-state operator}
\gianni{
\Beq
  \vbox{\hsize.7\hsize\noindent
  $\calS:\UR \to \calY , \quad
  u\mapsto\calS(u)=(\phi,\mu)$,
  \quad where \,$(\phi,\mu,\beta(\phi))\,$ 
  is the solution to the state
  system \Pbl\ corresponding to $u$,}
  \label{defS}
\Eeq
}%
{is well defined on $\UR$.}  

Based on \Start,
we prove in Theorem~\ref{Frechet} that $\calS$ is \Frechet\ differentiable at every point $\ustar\in\UR$,
where the \Frechet\ derivative is related 
to the solution to a linear system obtained from \Pbl\ by linearization.
{A standard argument then yields a necessary condition for optimal controls $\ustar\in\Uad$ that involves
the solution to the linearized system (Proposition~\ref{Badoptimality}).
As usual, we eliminate the solution to the linearized system by
means of the adjoint state variables $(p,q)$, 
which satisfy}
\begin{equation}\label{regpq}
  p \in \H1\Vp \cap \L\infty V \cap \L2W
  \aand
  q \in \L2V
\end{equation}
{and solve the adjoint state system}
	\begin{eqnarray}
	  && - \< \dt p(t) , v >
  + \iO p(t) {v}
  + \iO \nabla {q(t)} \cdot \nabla v
  + \iO f''(\phistar(t)) q(t) v
  \non
  \\
  && = \alpha_1 \iO (\phistar-\phiQ)(t) v
  \non
  \\
  && 
	\quad \hbox{\aat\ and every $v\in V$}\,,
  \label{primaA}
  \\[2mm]
  && \iO q(t) v
  = \iO \nabla p(t) \cdot \nabla v
  - \alpha_3 \iO (\mustar-\muQ)(t) v 
  \non
  \\
  && 
	\quad \hbox{\aat\ and every $v\in V$}\,,
  \label{secondaA}
  \\[2mm]
  && p(T)
  = \alpha_2 \bigl( \phistar(T) - \phiO \bigr)\,,
  \label{cauchyA}
\end{eqnarray}
\Accorpa\Adjoint primaA cauchyA
{where $(\phistar,\mustar)={\cal S}(\ustar)$.}
Our final result is the following necessary condition for optimality:

\Bthm
\label{Optimality}
Suppose that the conditions \Start\ are fulfilled, let $\ustar\in\Uad$ be an optimal
control with
associated state $(\phistar,\mustar)={\cal S} (u^*)$, and assume that
\begin{align}
&\alpha_1\phi_Q\in L^2(0,T;V),\quad \alpha_2\phi_\Omega\in V,\quad\alpha_3\mu_Q\in L^2(0,T;V),
\label{ende}\\
&\nabla\phistar\in L^\infty(0,T;L^4(\Omega)).
\label{finale}
\end{align}
Then it holds the variational inequality
\Beq
  \intQ p (u-\ustar) + \alpha_4 \intQ (u-\ustar)
  \geq 0
  \quad \forall\, u\in\Uad\,,
  \label{optimality}
\Eeq
where $p$ is the first component of the solution $(p,q)$ to the adjoint problem~\Adjoint.
\Ethm

\Brem
\gianni{Notice that each of the regularity conditions \eqref{ende} is fulfilled 
if either the related $\alpha_i\,$ is equal to $0$ 
or if the associated datum belongs to $L^2(0,T;V)$ or~$V$}.
In the first case, we have no corresponding tracking of the tumor fraction
$\phi$ or the chemical potential $\mu$, which is inconvenient from 
the medical viewpoint, while in the second case we assume the same
regularity for the respective target function that the optimal 
state $(\phistar,\mustar)$ is known to have
\gianni{(Theorem~\ref{Regularity})}, 
which seems to be reasonable.
\Erem

\vspace{2mm}
Throughout the paper, we will repeatedly use Young's inequality
\Beq
  a\,b \leq \delta\,a^2 + \frac 1{4\delta} \, b^2
  \quad \hbox{for all $a,b\in\erre$ and $\delta>0$},
  \label{young}
\Eeq
as well as \Holder's inequality and the Sobolev inequality 
related to the continuous embedding $V\subset L^p(\Omega)$ for $p\in[1,6]$
(since $\Omega$ is $d$-dimensional with $d\leq3$, bounded and smooth).
Furthermore, the embeddings $V\subset H$ and $H\subset\Vp$ are compact,
so that we {obtain from Ehrling's lemma} the compactness inequality 
\Beq
  {\|v\|}
  \,\leq\, \delta\, \norma{\nabla v}
  + C_\delta \, \norma v_{\Vp}
  \quad \hbox{for every $v\in V$ and $\delta>0$},
  \label{compact}
\Eeq
with some $C_\delta>0 $ that depends only on~$\Omega$ and~$\delta$.

We also utilize a tool that is widely used
in the study of Cahn--Hilliard type equations.
We define
\Beq
  \dom\calN := \graffe{\vstar\in\Vp: \ \barvstar = 0}
  \aand
  \calN : \dom\calN \to \graffe{v \in V : \ \barv = 0}
  \label{predefN}
\Eeq
by setting, for every $\vstar\in\dom\calN$,
\begin{align}
  & \hbox{$\calN\vstar$ is the unique element of $V$ that satisfies}
  \non
  \\
  & \iO \nabla\calN\vstar \cdot \nabla v = \< \vstar , v >
  \quad \hbox{for every $v\in V$}
  \aand 
  \overline{\calN\vstar} = 0 \,.
  \label{defN}
\end{align}
As $\Omega$ is bounded, smooth, and connected,
it turns out that \eqref{defN} yields a well-defined isomorphism,
which satisfies
\Beq
  \< \ustar , \calN \vstar >
  = \< \vstar , \calN \ustar >
  = \iO (\nabla\calN\ustar) \cdot (\nabla\calN\vstar)
  \quad \hbox{for every $\ustar,\vstar\in\dom\calN$}.
  \label{simmN}
\Eeq
Moreover, by also accounting for Poincar\'e's inequality
\Beq
  \normaV v^2 \leq \CO \bigl( \norma{\nabla v}^2 + |\barv|^2 \bigr)
  \quad \hbox{for every $v\in V$}\,,
  \label{poincare}
\Eeq
where the constant $\CO>0$ depends only on~$\Omega$,
we see that the function $\normaVp\cpto:\Vp\to[0,+\infty)$ defined by the formula
\Beq
  \normaVp\vstar^2
  := \norma{\nabla\calN(\vstar-\barvstar)}^2
  + |\barvstar|^2
  = \< \vstar , \calN(\vstar-\barvstar) >
  + |\barvstar|^2
  \quad \hbox{for $\vstar\in\Vp$},
  \label{defnormaVp}
\Eeq
is a Hilbert norm on $\Vp$ that is equivalent to the usual dual norm.
It follows that
\Beq
  | \< \vstar , v >|
  \leq \CO \normaVp\vstar \normaV v
  \quad \hbox{for every $\vstar\in\Vp$ and $v\in V$},
  \label{pier1}
\Eeq
with the same $\CO$ as in \eqref{poincare} without loss of generality.
Finally, notice that
\Beq
  2 \,\< \dt\vstar(t) , \calN\vstar(t) >
  = \frac d{dt} \iO |\nabla\calN\vstar(t)|^2
  = \frac d{dt} \, \normaVp{\vstar(t)}^2
  \quad \aat,
  \label{dtcalN}
\Eeq
\Accorpa\PropN predefN dtcalN
for every $\vstar\in\H1\Vp$ satisfying $\barvstar=0$, i.e.,
$\overline{\vstar(t)}=0$ \aat.
This kind of notation is used throughout the paper for time-dependent functions.

We conclude this section by stating a general rule concerning the constants 
that appear in the estimates to be performed in the following.
The small-case symbol $\,c\,$ stands for a generic constant
whose actual values may change from line to line and even within the same line
and depend only on~$\Omega$, on the shape of the nonlinearities,
and on the constants and the norms of the functions involved in the assumptions of the statements.
In particular, the values of $\,c\,$ do not depend on the parameters $\,\eps\,$ and $\,n\,$ 
we introduce in the next sections.
A~small-case symbol with a subscript like $c_\delta$ (in~particular, with $\delta=\eps$)
indicates that the constant may depend on the parameter~$\delta$, in addition.
On the contrary, we mark precise constants that we can refer~to
by using different symbols, e.g., capital letters like in~\eqref{compact}.


\section{Continuous dependence and uniqueness}
\label{CONTDEP}
\setcounter{equation}{0}

In this section, we give the proof of second part of Theorem~\ref{Wellposedness} 
and of Theorem~\ref{Lipcontdep}.
As for the former, we observe that the uniqueness 
of the \gianni{the component $\phi$ of the solution} follows from~\eqref{contdep1}
provided that this inequality is proved for every solution $(\phi_i,\mu_i,\xi_i)$
corresponding to~$u_i$ for $i=1,2$, \gianni{according to the statement}.
Indeed, \eqref{contdep1} with $u_1=u_2$ implies that $\phi_1=\phi_2$.
\gianni{Moreover, if $\beta$ is single-valued, \nju{then} we also deduce that $\xi_1=\xi_2$,
and a comparison in \eqref{seconda}, written for both solutions, yields that $\mu_1=\mu_2$ as well.}

\subsection{Continuous dependence, part I}
\label{CONTDEP1}

So we just have to prove the inequality \eqref{contdep1} by assuming that
$(\phi_i,\mu_i,\xi_i)$ is any solution corresponding to~$u_i$, $i=1,2$.
We set, for convenience, 
\Beq
  \phi := \phi_1-\phi_2 \,, \quad
  \mu := \mu_1-\mu_2 \,, \quad
  \xi := \xi_1-\xi_2
  \aand
  u := u_1-u_2 \,.
  \label{differenze}
\Eeq
Thus, by noting that \eqref{prima} and \eqref{ode} hold true with the new notations,
we can subtract the latter integrated over $\Omega$ from the former
and test the resulting equality by $\calN(\phi-\barphi)$
(recall \PropN\ for the definition of $\calN$ and its properties).
At~the same time, we write \eqref{seconda} for both solutions
and test the difference by $-(\phi-\barphi)$.
Then, we sum up and integrate over~$(0,t)$.
Since a cancellation occurs, we obtain~that
\Bsist
  && \frac 12 \, \normaVp{(\phi-\barphi)(t)}^2
  + \iot \normaVp{(\phi-
	\barphi)(s)}^2 \, ds
  + \intQt |\nabla(\phi-\barphi)|^2
  + \intQt \xi \, \phi
  \non
  \\
  && = \intQt (u-\baru) \, \calN(\phi-\barphi)
  - \intQt \xi \, \barphi
  - \intQt \bigl( \pi(\phi_1)-\pi(\phi_2)) (\phi-\barphi).
  \non
\Esist
All of the terms on the \lhs\ are nonnegative (the last one by monotonicity).
As for the first term on the \rhs, we account for \eqref{soluzode}
and notice that
$\norma\barphi_\infty\leq\norma\baru_{L^1(0,T)}\leq c\,\norma u_{\L1\Vp}$.
Hence, \pier{with the help of \eqref{pier1}} we have~that
\Bsist
  && \intQt (u-\baru) \, \calN(\phi-\barphi)
  \leq \pier{{}C_\Omega{}} \iot \normaVp{(u-\baru)(s)} \, \normaV{\calN(\phi-\barphi)(s)} \, ds
  \non
  \\
  && \leq c \iot \normaVp{(u-\baru)(s)} \, \normaVp{(\phi-\barphi)(s)} \, ds
  \leq c \, \norma u_{\L1\Vp} \sup_{0\leq s\leq t} \normaVp{(\phi-\barphi)(s)} 
  \non
  \\
  && \leq \frac 14 \sup_{0\leq s\leq t} \normaVp{(\phi-\barphi)(s)}^2
  + c \, \norma u_{\L1\Vp}^2 \,.
  \non
\Esist
Next, we have that
\Beq
  \intQt \xi \, \barphi
  \leq \bigl( \norma{\xi_1}_1 + \norma{\xi_2}_1 \bigr) \, \norma\barphi_\infty
  \leq c \bigl( \norma{\xi_1}_1 + \norma{\xi_2}_1 \bigr) \, \norma\baru_1 \,.
  \non
\Eeq
Finally, we owe to the \Lip\ continuity of $\pi$ and the compactness inequality \eqref{compact} to obtain~that
\Bsist
  && - \intQt \bigl( \pi(\phi_1)-\pi(\phi_2)) (\phi-\barphi)
  \leq c \intQt |\phi| \, |\phi-\barphi|  
  \non
  \\
  && \leq c \intQt |\phi-\barphi|^2
  + c \intQt |\barphi| \, |\phi-\barphi|
  \leq c \intQt |\phi-\barphi|^2
  + \iot |\barphi(s)|^2 \, ds
  \non
  \\
  && \leq \frac 12 \intQt |\nabla(\phi-\barphi)|^2
  + c \iot \normaVp{(\phi-\barphi)(s)}^2 \, ds
  + c \, \norma u_{\L1\Vp}^2\,.
  \non
\Esist
By collecting all these estimates, and ignoring some nonnegative terms on the \lhs, we deduce that
\Bsist
  && \frac 12 \, \normaVp{(\phi-\barphi)(t)}^2
  + \frac 12 \intQt |\nabla(\phi-\barphi)|^2
  \non
  \\
  && \leq \frac 14 \sup_{0\leq s\leq t} \normaVp{(\phi-\barphi)(s)}^2
  + c \iot \normaVp{(\phi-\barphi)(s)}^2 \, ds
  \non
  \\[2mm]
  && \quad {} 
  + c \, \norma u_{\L1\Vp}^2 
  + c \bigl( \norma{\xi_1}_1 + \norma{\xi_2}_1 \bigr) \, \norma\baru_1 \,.
  \non
\Esist
By applying Gronwall's lemma, we conclude that 
\eqref{contdep1} holds true with a constant $C_1$ \nju{as} in the statement.

\subsection{Continuous dependence, part II}

We are going to prove Theorem~\ref{Lipcontdep}.
We still use the notations \eqref{differenze} regarding $\phi$, $\mu$, and~$u$,
and denote the \Lip\ constant of $f'$ by~$L$.
We write \eqref{prima} and \eqref{seconda} for both solutions
and test the differences by $\phi$ and~$\mu$, respectively.
Then, \pier{when adding} an obvious cancellation occurs, and 
\pier{by Young's inequality} we have that
\Bsist
  && \frac 12 \iO |\phi(t)|^2
  + \intQt |\phi|^2
  + \intQt |\mu|^2
  \non
  \\
  && = \intQt u \phi
  + \pier{{}\intQt{}} \bigl( f'(\phi_1) - f'(\phi_2) \bigr) \pier{{}\mu{}}
  \non
  \\
  &&\leq \pier{\frac12} \intQ |u|^2
  + \pier{\frac{1+L^2}2 \intQt |\phi|^2 + \frac{1}2 \intQt |\mu|^2}\,.
  \non
\Esist
\pier{Hence, rearranging and using} Gronwall's lemma, we immediately deduce that
\Beq
  \norma\phi_{\C0H} 
  + \norma\mu_{\L2H}
  \leq c \, \norma u_{\L2H}\,,
  \non
\Eeq
with a constant $c$ that only depends on $L$ and~$T$.
By applying elliptic regularity to the difference of \eqref{seconda}, written for both solutions,
we also infer~that
\Bsist
  && \nju{\|\phi\|_{L^2(0,T;W)}}
  \leq c \, \norma{\mu-(f'(\phi_1)-f'(\phi_2))}_{\L2H}
  \non
  \\[2mm]
  && \leq c \bigl( \norma\mu_{\L2H} + \norma\phi_{\L2H} \bigr)
  \leq c \, \norma u_{\L2H}\,,
  \non
\Esist
where $c$ depends on~$\Omega$, in addition.


\section{Approximation and existence}
\label{EXISTENCE}
\setcounter{equation}{0}

In this section, we prove the existence part of Theorem~\ref{Wellposedness}.
Our method consists in approximating the problem at hand and using compactness and monotonicity arguments.

\subsection{Approximation}
Here, we construct an approximating problem
depending \gianni{on the parameter $\eps\in(0,1)$.
This problem is obtained by modifying the formulation of problem \Pbl\
by replacing $\Beta$ and $\beta$ by their Moreau-Yosida regularizations
$\Betaeps$ and~$\betaeps$, respectively}
(see, e.g., \cite[pp.~28 and~39]{Brezis}).
\pier{Owing also} to our assumption \eqref{hpBeta} on~$\Beta$, we have~that
\Bsist
  && \hbox{$\betaeps$ is monotone and \Lip\ continuous with $\betaeps(0)=0$},
  \label{monbetaeps}
  \\
  && |\betaeps(r)| \leq |\betaz(r)|
  \quad \hbox{for every $r\in D(\beta)$},
  \label{disugbetaeps}
  \\
  && 0 \leq \Betaeps(r) = \int_0^r \betaeps(s) \, ds \leq \Beta(r)
  \quad \hbox{for every $r\in\erre$},
  \label{disugBetaeps}
\Esist
\Accorpa\Propbetaeps monbetaeps disugBetaeps
\gianni{where we recall that $\betaz(r)$ denotes the element of $\beta(r)$ having minimal modulus.}
The solution we look for is a pair $(\phieps,\mueps)$
enjoying the regularity properties 
\Beq
  \gianni{\phieps \in \H1\Vp \cap \L\infty V \pier{,} \quad 
  \mueps \in \L2V}
  \label{regsoluzeps}
\Eeq
and satisfying 
\Bsist
  && \iO \dt\phieps(t) \, v
  + \iO \phieps(t) v 
  + \iO \nabla\mueps(t) \cdot \nabla v
  = \iO u(t) v
  \non
  \\
  && \quad \hbox{\aat\ and every $v\in V$},
  \label{primaeps}
  \\
  && \iO \mueps(t) v 
  = \iO \nabla\phieps(t) \cdot \nabla v
  + \iO \betaeps(\phieps(t)) v
  + \iO \pi(\phieps(t)) v
  \qquad
  \non
  \\
  && \quad \hbox{\aat\ and every $v\in V$},
  \label{secondaeps}
  \\
  && \phieps(0) = \phiz \,.
  \label{cauchyeps}
\Esist
\Accorpa\Pbleps primaeps cauchyeps
Our aim is \nju{to solve} this approximating problem.
\nju{In} this respect, we have the following result.

\Bthm
\label{Wellposednesseps}
Assume \HPstruttura\ on the structure on the system and \HPdati\ on the data.
Then, for every \gianni{$\eps\in(0,1)$}, 
problem \Pbleps\ has a unique solution $(\phieps,\mueps)$
satisfying the regularity properties~\eqref{regsoluzeps}.
\Ethm

\gianni{Uniqueness follows from Section~\ref{CONTDEP1},
since $\betaeps$ satisfies all the properties required for~$\beta$
and is single-valued}.
Hence, only the existence of a solution to problem \Pbleps\ has to be proved.
Our argument is based on the discretization of \Pbleps\ by a Faedo--Galerkin scheme
and a~priori estimates.

\step
The discrete problem

We introduce the nondecreasing sequence $\{\lambda_j\}$ of the eigenvalues
and the complete orthonormal sequence $\{e_j\}$ of corresponding eigenvectors
of the Laplace operator with homogeneous Neumann boundary conditions,
that~is,
\begin{align}
  &\juerg{-\Delta e_j = \lambda_j e_j
  \quad \hbox{in $\Omega$}
  \aand
  \dn e_j = 0 \quad\mbox{on \,$\Gamma$}\quad \mbox{for all \,$j\in\enne$},
	\quad\mbox{with} \quad
  (e_i , e_j) = \delta_{ij}}\non\\
	&\,\juerg{\mbox{and}\quad (\nabla e_i,\nabla e_j)=\lambda_i\,\delta_{ij}
  \quad \hbox{for all}\quad i,j\in\enne.}
  \label{eigen}
\end{align}
Next, for $n\geq1$, we introduce the subspace $\Vn$ of $V$ by setting
\Beq
  \Vn := \Span \graffe{e_1,\dots,e_n}.
  \label{defVn}
\Eeq
Then, the sequence $\graffe{\Vn}$ is nondecreasing, and its union is dense in both $V$ and~$H$.
We observe at once that 
\Beq
  \hbox{$\lambda_1=0$ and the elements of $V_1$ are the constant functions}.
  \label{propV1}
\Eeq
At this point, we can introduce the discrete problem.
Even though it also depends \gianni{on~$\eps$}, 
we do not stress this in the notation for the solution.
We look for a pair $(\phin,\mun)$ satisfying
\Bsist
  && \gianni{\phin,\, \mun \in \H1\Vn},
  \label{regsoluzn}
  \\[2mm]
  && \juerg{(\dt\phin(t), v) 
  + (\phin(t), v) 
  + (\nabla\mun(t),\nabla v)
  = (u(t), v)}
  \non
  \\
  &&\juerg{\quad \hbox{\aat\ and every $v\in\Vn$}},
  \label{priman}
  \\[1mm]
  && \gianni{(\mun(t), v) 
  = ( \nabla\phin(t), \nabla v)
  + (\betaeps(\phin(t)), v)
  + (\pi(\phin(t)) ,v)}
  \non
  \qquad
  \\
  && \juerg{\quad \hbox{for every $t\in[0,T]$ and $v\in\Vn$}},
  \label{secondan}
  \\[1mm]
  && \juerg{(\phin(0), v)
  = (\phiz, v) \quad \hbox{for every $v\in\Vn$}.}
  \label{cauchyn}
\Esist
\Accorpa\Pbln priman cauchyn

\Brem
\label{Phinz}
We notice that $\phin(0)$ is the $H$-projection of $\phiz$ on $\Vn$, i.e.,
\Beq
  \phin(0) = \somma j1n (\phiz,e_j) \, e_j
  \quad \hbox{for $n=1,2,\dots$}.
  \label{phinz}
\Eeq 
This and assumption \eqref{hpphiz} have important consequences.
The first one is given by the inequalities
\Bsist
  && \norma{\phin(0)} \leq \norma\phiz
  \aand
  \non
  \\
  && \quad 
  \norma{\nabla\phin(0)}^2 
  \leq c \somma j1n |\lambda_j^{1/2}(\phiz,e_j)|^2
  \leq c \somma j1\infty |\lambda_j^{1/2}(\phiz,e_j)|^2
  \leq c \, \normaV\phiz^2 \,.
  \label{normephinz}
\Esist
Next, since $\phiz\in W$, we have that
\begin{align}
  &\somma j1\infty |\lambda_j(\phiz,e_j)|^2
  \leq c \, \normaW\phiz^2
  \aand\non\\
  &\normaW{\phiz-\phin(0)}^2
  \leq c \somma j{n+1}\infty |\lambda_j(\phiz,e_j)|^2
  \nju{\to 0\quad\mbox{as $n\to\infty$}}, 
  \non
\end{align}
so that $\phin(0)$ converges to $\phiz$ in~$W$, and thus uniformly.
By combining with our assumption \eqref{hpphiz} on $\phiz$, 
we deduce that there are elements $r'_\pm$ of the interior of $D(\beta)$
and a natural number $n_0$ only depending on $\beta$ and $\phiz$ such that
$r'_-\leq\phin(0)\leq r'_+$ in $\,\Omega\,$ for every $n\geq n_0$.
By accounting for \accorpa{disugbetaeps}{disugBetaeps}, we infer that
\Beq
  \Betaeps(\phin(0)) + |\betaeps(\phin(0)\pier)|
  \leq \sup_{r'_-\leq r\leq r'_+} \Beta(r)
  + \sup_{r'_-\leq r\leq r'_+} |\betaz(r)| 
  \label{betaepsphinz}
\Eeq
for every $n\geq n_0$.
From now on, it is understood that $n\geq n_0$ so that \eqref{betaepsphinz} holds true.
\Erem

The next step is proving that this problem has a unique solution.
To this end, we represent $\phin$ and $\mun$ in term of the eigenfunctions, i.e.,
\Beq
  \phin(t) = \somma j1n \phinj(t) \, e_j
  \aand
  \mun(t) = \somma j1n \munj(t) \, e_j
  \quad \hbox{for every $t\in[0,T]$},
  \non
\Eeq
\gianni{for some functions $\phinj$ and $\munj$ belonging to $H^1(0,T)$}.
Then, if we consider the column vectors
$y:=(\phinj)_{j=1,\dots,n}$ and $z:=(\munj)_{j=1,\dots,n}$,
equations \accorpa{priman}{secondan} take the form
\Beq
  y'(t) + y(t) + A z(t) = g(t)
  \aand
  \gianni{z(t) = A y(t) + G(y(t))},
  \non
\Eeq
where $A$ is the diagonal matrix of the first $n$ eigenvalues,
$g(t)$ is the column vector with components \gianni{$(u(t),e_j)$},
and $G:\erren\to\erren$ is a \Lip\ continuous function.
Notice that $g\in\L2{\erre^n}$.
Hence, by replacing $z$ in the first equation with the help of the second one and rearranging,
since \eqref{phinz} provides an initial condition for~$y$,
we obtain a well-posed Cauchy problem for~$y$ which has a unique solution $y\in\H1{\erre^n}$.
\gianni{Then, the second equation yields that $z\in\H1{\erre^n}$.}

\subsection{Uniform a priori estimates}
\label{UNIFEST}

We perform a number of estimates.
We point out that they are based on the properties \Propbetaeps\ 
and thus hold if $\betaeps$ is a possibly different approximation of~$\beta$.
We make this remark since the estimates we are going to prove here 
will be used later on for other purposes with a different~$\betaeps$.
Hence, we recall our general rule:
the symbol $c$ denotes (possibly different) constants that do not depend on \gianni{$\eps$ and~$n$}.

However, before starting, we remark a consequence of~\eqref{compatibility}.
We choose $\delta_0>0$ such that 
the interval $[\barphiz-M-\delta_0,\barphiz+M+\delta_0]$ 
is included in the interior of~$D(\beta)$.
Then, for some $C_0>0$, we have the inequality
\Bsist
  && \betaeps(r) (r-r_0)
  \geq \delta_0 |\betaeps(r)| - C_0
  \non
  \\
  && \quad \hbox{for every $r\in\erre$, $r_0\in[\barphiz-M,\barphiz+M]$ and $\eps\in(0,1)$}.
  \label{trickMZ}
\Esist
This is a generalization of \cite[Appendix, Prop. A.1]{MiZe}.
The detailed proof given in \cite[p.~908]{GiMiSchi} with a fixed $r_0$ 
also works in the present case with minor changes.

\step
Preliminary estimate

We notice that, by \eqref{propV1}, the constant functions belong to $\Vn$ for every~$n$.
Hence, we deduce from \eqref{priman} and \eqref{cauchyn} that
\Beq
  \frac d{dt} \, \barphin(t) + \barphin(t) = \baru(t)
  \quad \aat
  \aand
  \barphin(0) = \barphiz \,.
  \label{oden}
\Eeq
As in Remark~\ref{Compatibility} (both the equation and the initial condition are the same here and there), we conclude that
\Beq
  \barphiz - M \leq \barphin(t) \leq \barphiz + M
  \quad \hbox{for every $t\in[0,T]$} .
  \label{stimabarphin}
\Eeq

\step
First uniform estimate

We recall the definition of $\calN$ and its properties (see \PropN) and notice that
\Beq
  v \in \Vn \aand \barv =0 
  \quad \hbox{imply that} \quad
  \calN v \in \Vn \,.
  \label{pier2}
\Eeq
Indeed, both $v$ and $w:=\calN v$ can be expressed in terms of the eigenfunctions~$e_j$, 
and we have that 
\Beq
  \somma n1\infty \lambda_j (w,e_j) \, e_j 
  = -\Delta w
  = v
  = \somma j2n (v,e_j) \, e_j \,.
  \non
\Eeq
Hence, $(w,\nju{e_j})=0$ for every $j>n$ (since $\lambda_j>0$ for $j>1$), i.e.,  $w\in\Vn$.
Therefore, we subtract \eqref{oden}, multiplied by $\iO v$, \nju{from} \eqref{priman}, 
write the resulting equality at the time $s$ 
and test it by $\calN(\phin(s)-\barphin(s))$.
At the same time, we test \eqref{secondan} written at the time $s$ by $-(\phin(s)-\barphin(s))$.
Then, we sum up and integrate over $(0,t)$ with respect to~$s$.
By accounting for the cancellation that occurs, we obtain~that
\gianni{%
\Bsist
  && \frac 12 \, \normaVp{\phin(t)-\barphin(t)}^2
  + \iot \normaVp{\phin(s)-\barphin(s)}^2 \, ds
  \non
  \\
  && \quad {}
  + \intQt |\nabla\phin|^2
  + \intQt \betaeps(\phin)(\phin-\barphin)
  \non
  \\
  && = \frac 12 \, \normaVp{\phin(0)-\barphiz}^2
  + \intQt (u-\baru)\,\pier{\calN}(\phin-\barphin)
  - \intQt \pi(\phin)(\phin-\barphin) \,.
  \non
\Esist
}%
The first three terms on the \lhs\ are nonnegative.
To treat the last one, we account for \eqref{trickMZ} and \eqref{stimabarphin} and have that
\Beq
  \iO \betaeps(\phin(t)) (\phin(t)-\barphin(t))
  \geq \delta_0 \, \norma{\betaeps(\phin(t))}_{\Luno} - c 
  \quad \aat , 
  \label{fromMZ}
\Eeq
whence also
\Beq
  \intQt \betaeps(\phin)(\phin-\barphin)
  \geq \delta_0 \norma{\betaeps(\phin)}_{\pier{L^1(Q_t)}} - c \,.
  \non
\Eeq
Let us come to the \rhs.
The first \gianni{term is} estimated from above by using the $H$ norm.
Thus, \gianni{it is} uniformly bounded since $\pier{\phi_n}(0)$ is the $H$-projection of $\phiz$ on~$\Vn$.
\gianni{The next term can be dealt with} \pier{as follows:}
\begin{align}
  &\intQt \gianni{(u-\baru)}\,\pier{\calN}(\phin-\barphin)
  \pier{{}\leq  \iot \normaVp{u(s)-\baru(s)} \normaVp{\phin(s)-\barphin(s)} ds  }
  \non \\ &
  \leq  \pier{c \norma{u-\baru}_{\L2H}^2 
  + \frac12  \iot \normaVp{\phin(s)-\barphin(s)}^2 ds\,.}
  \non
\end{align}
By adding and subtracting $\pi(\barphin)$ in the first factor inside the last integral,
and accounting for the \Lip\ continuity of~$\pi$,
we have that
\Beq
  - \intQt \pi(\phin)(\phin-\barphin)
  \leq c \nju{\intQt} |\phin-\barphin|^2
  + c \intQt |\barphin|^2 + c \,.
  \non
\Eeq
So, we can treat the first integral on the \rhs\ 
\pier{with the help of the compactness inequality~\eqref{compact}, namely 
\Beq
  c {\intQt} |\phin-\barphin|^2
  \leq  \delta \intQt |\nabla\phin|^2
  + c_\delta \iot \normaVp{\phin(s)-\barphin(s)}^2 ds\,,
  \non
\Eeq
where $\delta>0$ is arbitrary,}
and we can use \eqref{stimabarphin} for the second one.
At this point, we choose $\delta>0$ small enough and apply the Gronwall lemma.
We conclude~that
\gianni{%
\Beq
  \norma\phin_{\L\infty\Vp\cap\L2V}
  + \norma{\betaeps(\phin)}_{\LQ1}
  \leq c \,.
  \label{primastiman}
\Eeq
}%
Since $\barmun=\overline{\betaeps(\phin)}+\overline{\pi(\phin)}$
by \eqref{secondan} tested  by $1/|\Omega|$,
we deduce that
\Beq
  \norma\barmun_{L^1(0,T)} \leq c \,.
  \label{daprimastiman}
\Eeq

\step
Second uniform estimate

We test \eqref{priman} by $\mun$ and \eqref{secondan} by $-(\dt\phin+\phin)$ and sum up.
Then, an obvious \nju{cancellation} occurs, and \pier{integration over $(0,t)$ leads to}
\gianni{%
\Bsist
  && \intQt |\nabla\mun|^2
  + \frac 12 \iO |\nabla\phin(t)|^2
  + \intQt |\nabla\phin|^2
  + \iO \Betaeps(\phin(t))
  + \intQt \betaeps(\phin) \phin
  \non
  \\
  && = \frac 12 \iO |\nabla\phin(0)|^2
  + \iO \Betaeps(\phin(0))
  \non
  \\
  && \quad {}
  + \intQt u \mun
  - \iO \Pi(\phin(t))
  + \iO \Pi (\phin(0))
  - \intQt \pi(\phin) \phin \,.
  \non
\Esist
}%
\nju{All of} the terms on the \lhs\ are nonnegative.
The first \gianni{two} terms on the \rhs\ can be estimated \pier{thanks} to \accorpa{normephinz}{betaepsphinz}.
The same \eqref{normephinz}\pier{, along with \eqref{hpPi},} helps in the second term involving~$\Pi$,
while the last integral is bounded due to \eqref{primastiman}.
It remains to treat two \nju{further} terms.
By accounting for Young's inequality, the Poincar\'e inequality \eqref{poincare}, and \eqref{daprimastiman},
we have~that
\Bsist
  && \intQt u \mun
  = \intQt u (\mun-\barmun)
  + \intQt u \barmun
  \non
  \\
  && \leq \frac 12 \intQt |\nabla\mun|^2 
  + c \intQt |u|^2
  + |\Omega| \, \norma u_\infty \, \norma\barmun_{L^1(0,T)}
  \leq \frac 12 \intQt |\nabla\mun|^2 
  + c \,.
  \non
\Esist
Finally, we owe to the quadratic growth of~$\Pi$, 
the compactness inequality \eqref{compact}, and \eqref{primastiman},
to see that
\Beq
  - \iO \Pi(\phin(t))
  \leq c \iO |\phin(t)|^2
  + c
  \leq \frac 14 \iO |\nabla\phin(t)|^2
  + c \, \normaVp{\phin(t)}^2 
  + c 
  \leq \frac 14 \iO |\nabla\phin(t)|^2
  + c \,. 
  \non
\Eeq
By combining all these estimates, we conclude that
\Beq
  \gianni{\norma{\nabla\mun}_{\L2H}
  + \norma\phin_{\L\infty V}
  + \norma{\Betaeps(\phin\pier{)}}_{\L\infty\Luno}
  \leq c \,.}
  \label{secondastiman}
\Eeq

\step
Third uniform estimate

We test \eqref{secondan} by $\phin(t)-\barphin(t)$.
We have \aet\ that
\Bsist
  && \iO |\nabla\phin|^2
  + \iO \betaeps(\phin)(\phin-\barphin)
  \non
  \\
  && = - \iO \pi(\phin) (\phin-\barphin)
  + \iO \gianni\mun (\phin-\barphin) .
  \non
\Esist
The second term on the \lhs\ is estimated from below by \eqref{fromMZ}, \pier{and} 
the first term on the \rhs\ is \pier{uniformly bounded by \eqref{secondastiman}
and the Lipschitz continuity of $\pi$. The
last term} is treated with the help of Poincar\'e's inequality:
\Bsist
  && \iO \gianni\mun (\phin-\barphin)
  = \iO \gianni{(\mun-\barmun)} (\phin-\barphin)
  \non
  \\
  && \leq c \, \norma{\phin-\barphin} \gianni{\, \norma{\nabla\mun}}
  \leq c \gianni{\, \norma{\nabla\mun}} \,.
  \non
\Esist
By combining, we deduce that 
\Beq
  \delta_0 \norma{\betaeps(\phin(t))}_{\Luno}
  \leq c \gianni{\, \norma{\nabla\mun(t)}}
  + c \,,
  \quad \hbox{i.e.,} \quad
  \delta_0^2 \norma{\betaeps(\phin(t))}_{\Luno}^2
  \leq c \gianni{\, \norma{\nabla\mun(t)}^2}
  + c
  \non
\Eeq
\aat, and \eqref{secondastiman} yields that
\gianni{%
\Beq
  \norma{\betaeps(\phin)}_{\L2\Luno}
  \leq c \, \norma{\nabla\mun}_{\L2H}
  + c 
  \leq c \,.
  \non
\Eeq
}%
Consequently, at first the mean value $\overline{\betaeps(\phin)}$, 
and then $\barmun$, are bounded in $L^2(0,T)$.
By accounting for \eqref{secondastiman} \pier{and the Poincar\'e inequality \eqref{poincare}, we then} conclude that
\Beq
  \norma\mun_{\L2V} \leq c \,.
  \label{terzastiman}
\Eeq

\subsection{Proof of Theorem~\ref{Wellposednesseps}}
\label{ProofWPeps}

The above estimates provide some weak convergence for both $\phin$ and~$\mun$.
However, \gianni{no strong convergence can \nju{as yet be deduced}}.
Thus, we perform one more a priori estimate, which is uniform at least with respect to~$n$.

\step
An auxiliary estimate

First, we subtract \eqref{oden}, multiplied by~$\iO v$, \nju{from} \eqref{priman}.
Then, we test the equality obtained by $\calN w$ and \eqref{secondan} by $-w$,
where $w:=\dt(\phin-\barphin)+(\phin-\barphin)$.
We integrate  the resulting equalities over $(0,t)$, rearrange the second one, and obtain that
\Bsist
  && \iot \normaVp{\dt(\phin(s)-\barphin(s))+(\phin(s)-\barphin(s))}^2 \, ds
  + \intQt \nabla\mun \cdot \nabla\calN w
  \non
  \\
  && = \intQt (u-\baru) \calN w\,,
  \non
  \\[2mm]
  \separa
  && \gianni{\iO |\nabla\phin(t)|^2 
  + \intQt |\nabla\phin|^2
  + \iO \Betaeps(\phin(t))}
  \non
  \\
  && = \gianni{\iO \Betaeps(\phin(0))
  - \iO \Pi(\phin(t))
  + \iO \Pi(\phin(0))}
  \non
  \\
  && \quad {}
  + \intQt \bigl( \betaeps(\phin)+\pi(\phin) \bigr) (-\phin+\dt\barphin+\barphin)
  + \intQt \mun w\,.
  \non
\Esist
At this point, we sum up and notice that a cancellation occurs due to the definition of~$\calN$.
Moreover, we owe to \eqref{oden} on both sides, and account for the quadratic growth of~$\Pi$,
the linear growth of $\betaeps+\pi$, and some of the inequalities \eqref{normephinz} and \eqref{betaepsphinz}.
We obtain that
\Bsist
  && \iot \normaVp{\dt\phin(s)+\phin(s)-\baru(s)}^2 \, ds
  + \iO |\nabla\phin(t)|^2
  + \intQt |\nabla\phin|^2
  + \iO \Betaeps(\phin(t))
  \non
  \\
  && \gianni{\leq \intQt (u-\baru) \calN w
  + c \iO |\phin(t)|^2
  + c_\eps \intQt \bigl( |\phin|^2 + |\baru|^2 \bigr)
  + c_\eps \,.}
  \non
\Esist
Thanks to \eqref{secondastiman} and the assumptions on~$u$, all of 
the terms on the \rhs\ are under control but the first one.
\nju{In} this respect, we have that
\Bsist
  && \intQt (u-\baru) \calN w
  \leq \norma{u-\baru}_{\L2H} \norma{\calN w}_{\Lt2H}
  \leq c \norma w_{\Lt2\Vp}
  \non
  \\
  && = c \, \norma{\dt\phin+\phin-\baru}_{\Lt2\Vp}
  \leq \frac 12 \iot \normaVp{\dt\phin(s)+\phin(s)-\baru(s)}^2 \, ds
  + c \,.
  \non
\Esist
By combining, rearranging, and noting that $\phin-\baru$ is obviously bounded in $\L2\Vp$,
we conclude that
\Beq
  \norma{\dt\phin}_{\L2\Vp} 
  \leq c_\eps \,.
  \label{auxil}
\Eeq

\step
Conclusion

Now, we can let $n$ tend to infinity in the discrete problem.
Indeed, the estimates \eqref{primastiman} and \accorpa{secondastiman}{terzastiman}
are uniform with respect to~$n$,
so that, by also applying well-known compactness results 
(for \nju{the strong convergence}, see, e.g., \cite[Sect.~8, Cor.~4]{Simon}),
we~\nju{conclude}  that there exists a pair $(\phieps,\mueps)$ such~that
\gianni{%
\Bsist
  & \phin \to \phieps
  & \quad \hbox{weakly star in $\gianni{\H1{V^*}}\cap\L\infty V$}
  \non
  \\
  && \qquad \hbox{and strongly in $\C0H$}\,,
  \label{convphin}
  \\
  & \mun \to \mueps
  & \quad \hbox{weakly in $\L2V$}\,,
  \label{convmun}
\Esist
}%
just for a subsequence, in principle.
However, since we prove that the limit $(\phieps,\mueps)$ is a solution to \Pbleps\
and we already know that uniqueness holds for this problem,
the whole sequences converge to $(\phieps,\mueps)$.
By \Lip\ continuity, we infer that
\Beq
  \betaeps(\phin) + \pi(\phin) \to \betaeps(\phieps) + \pi(\phieps)
  \quad \hbox{strongly in $\C0H$} .
  \non
\Eeq
Since $\phin(0)$ converges to both $\phieps(0)$ and $\phiz$ in~$H$,
the initial condition \eqref{cauchyeps} is satisfied,
and it remains to prove that the equations \accorpa{primaeps}{secondaeps} are satisfied as well.
It suffices to verify that $(\phieps,\mueps)$ solves the corresponding integrated versions
with time dependent test functions in $\L2V$,~i.e., 
\Bsist
  && \intQ \dt\phieps \, v
  + \intQ \phieps v
  + \intQ \nabla\mueps \cdot \nabla v
  = \intQ u v
  \non
  \\
  && \quad \hbox{for every $v\in\L2V$}\,,
  \label{intprimaeps}
  \\
  && \gianni{\intQ \mueps v 
  = \intQ \nabla\phieps \cdot \nabla v
  + \intQ \betaeps(\phieps) v
  + \iO \pi(\phieps) v}
  \non
  \\
  && \quad \hbox{for every $v\in\L2V$}\,.
  \label{intsecondaeps}
\Esist
\gianni{To this end, we fix $m$ for a while and take any $V_m$-valued step function~$w$.
Then, if $n\geq m$, the choice $v=w(t)$ is admissible in both \eqref{priman} and \eqref{secondan}.
By testing the equations and integrating over $(0,T)$, 
we obtain \accorpa{intprimaeps}{intsecondaeps} for $(\phin,\mun)$ with this test function,
and letting $n$ tend to infinity we have the same for $(\phieps,\mueps)$.
Since $w$ and $m$ are arbitrary, the equations \accorpa{intprimaeps}{intsecondaeps} also hold 
with any $V_\infty$-valued step function, where $V_\infty$ is the union of all the subspaces~$V_m$.
As $V_\infty$ is dense in~$V$, the set of these step functions is dense in $\L2V$, and we immediately conclude.}

\subsection{Proof of Theorem~\ref{Wellposedness}}
\label{ProofWP}

Since the estimates of Section~\ref{UNIFEST} are uniform with respect to \gianni{$n$ and~$\eps$},
they are preserved in the limit as $n$ tends to infinity, and we have that
\Beq
  \gianni{\norma\phieps_{\L\infty V}
  + \norma\mueps_{\L2V}
  + \norma{\Betaeps(\phieps)}_{\L\infty\Luno}
  \leq c \,.}
  \label{fromunifest}
\Eeq
\gianni{However, we need some \nju{further} estimates.}

\step
\pier{Fourth a priori estimate}

\pier{From \eqref{intprimaeps} it easily follows that 
\begin{align}
&\left| \intQ \dt\phieps \, v \right|
  \leq \Bigl(\norma\phieps_{\L\infty H}  + \norma{\nabla \mueps}_{\L2 H }  
  + \norma{u}_{\L2 H} \Bigr) \norma{v}_{\L2V}    
  \non \\
  &\quad \hbox{for every $v\in\L2V$},
\non
\end{align}
whence the estimate \eqref{fromunifest} implies that}
\Beq
  \norma{\dt\phieps}_{\L2\Vp} \leq c \,.
  \label{quintastima}
\Eeq

\step
\pier{Fifth a priori estimate}

We test \eqref{secondaeps} by $\betaeps(\phieps(t))$ and integrate over~$(0,T)$.
We obtain~that
\gianni{%
\Beq
  \intQ \betaeps'(\phieps) |\nabla\phieps|^2
  + \intQ |\betaeps(\phieps)|^2
  = \intQ \bigl( \mueps - \pi(\phieps) \bigr) \betaeps(\phieps) \,.
  \non
\Eeq
}%
By accounting for \eqref{disugBetaeps}, the \Lip\ continuity of $\pi$, and~\eqref{fromunifest},
we immediately conclude~that
\Beq
  \norma{\betaeps(\phieps)}_{\L2H} \leq c \,.
  \label{quartastima}
\Eeq

\step
Conclusion

\gianni{We \pier{deduce}~that
\Bsist
  & \phieps \to \phi
  & \quad \hbox{weakly star in $\H1\Vp\cap\L\infty V$}
  \non
  \\
  && \qquad \hbox{and strongly in $\C0H$}\,,
  \label{convphieps}
  \\
  & \mueps \to \mu
  & \quad \hbox{weakly in $\L2V$}\,,
  \label{convmueps}
  \\
  & \betaeps(\phieps) \to \xi
  & \quad \hbox{weakly in $\L2H$}\,,
  \label{convxieps}
\Esist
\Accorpa\Convsoluzeps convphieps convxieps
as $\eps$ tends to~$0$
(at~least for a subsequence $\eps_k$ tending to~$0$),
where $(\phi,\mu,\xi)$ is a solution to problem \Pbl.}
By virtue of the estimates \accorpa{fromunifest}{quintastima} and compactness results,
we see that \accorpa{convphieps}{convxieps} hold true
for some triplet $(\phi,\mu,\xi)$.
Then, \pier{in view of \eqref{cauchyeps}, it turns out that} 
\eqref{cauchy} is satisfied, and we can take the limit in 
\accorpa{intprimaeps}{intsecondaeps}
to obtain \accorpa{intprima}{intseconda}.
It remains to show that $\xi\in\beta(\phi)$ \aeQ,
but this follows from the strong convergence of $\phieps$,
the weak convergence of $\betaeps(\phieps)$, and a well-known property of the Yosida approximation
(see, e.g.,  \cite[Prop.~2.2, p.~38]{Barbu}).
Thus, the proof of Theorem~\ref{Wellposedness} is complete.


\section{Regularity and separation}
\label{REGULARITY}
\setcounter{equation}{0}

This section is devoted to the proofs of Theorem~\ref{Regularity}
and Propositions \ref{Separation} and~\ref{Separation2d}.
Hence, it is understood that the assumptions of the corresponding statements are in force.

\subsection{Regularity}

In order to prove Theorem~\ref{Regularity}, we make a preliminary remark.
Since $\beta$ is a $C^1$ function in the interior of its domain,
we can prove an estimate for the derivative $\betaeps'$ of the Yosida approximation.
By the definition of $\betaeps$, we have that
\Beq
  \betaeps(r) \in \beta(r-\eps\betaeps(r))
  \quad \hbox{for every $r\in\erre$}.
  \label{defbetaeps}
\Eeq
Assume now that $r$ belongs to the interior of~$D(\beta)$.
Then the same is true for $r-\eps\betaeps(r)$, and it holds that $|r-\eps\betaeps(r)|\leq|r|$.
Indeed, by supposing, e.g., that $r>0$, \nju{it} follows that $\betaeps(r)\geq0$ whence $r-\eps\betaeps(r)\leq r$.
\gianni{But we also have that $r-\eps\betaeps(r)\geq0$ 
since $\betaeps(0)=0$ and $\betaeps$ is \Lip\ continuous with \Lip\ constant~$1/\eps$.
This proves the second assertion.
The first one follows since $[0,r]$ is included in the interior of~$\beta$ (recall \eqref{hpbetareg}).}
Therefore, \eqref{defbetaeps} becomes an equation,
and $\pier{s=\betaeps(r)}$ is the \nju{function implicitly} defined by the equation $\beta(r-\eps s)-s=0$.
Since $(\partial/\partial s)(\beta(r-\eps s)-s)\leq-1$, $\betaeps$~is differentiable
and its derivative is given by 
\Beq
  \betaeps'(r) = \frac {\beta'(r-\eps\betaeps(r))} {1+\eps\beta'(r-\eps\betaeps(r))}\,.
  \non
\Eeq
Hence, we deduce that
\Beq
  |\betaeps'(r)| \leq \sup_{s\in[a,b]} |\beta'(s)|
  \quad \hbox{for every $r\in[a,b]$}\,,
  \label{stimabetaprimo}
\Eeq
whenever $[a,b]$ contains $0$ and is included in the interior of $D(\beta)$.

\step
Regularity estimate

To be completely rigorous, we come back to the \gianni{Faedo--Galerkin scheme}
and notice that \gianni{we can assume $\phin$ and $\mun$ to be of class~$\pier{C^1}$}.
Indeed, $u$~belongs to $\H1H$, and we can \gianni{suppose} that $\betaeps$ is of class~$C^1$
since it could be replaced by a regularization of it 
having the same properties in the whole procedure developed so far
(see, e.g., \cite[formulas (3.5--6)]{GR3}).
We denote the \Lip\ constant of $\pi$ by~$L$, set $L':=L+1$,
and test \eqref{priman} by $\dt\mun+L'\,\dt\phin$.
At the same time, we differentiate \eqref{secondan} 
and test the resulting equation by~$-(\dt\phin+\phin)$.
Then, we sum up, integrate over~$(0,t)$, and notice that a cancellation occurs.
Thus, we obtain that
\Bsist
  && \frac 12 \iO |\nabla\mun(t)|^2
  + L' \intQt |\dt\phin|^2
  + \frac {L'}2 \iO |\phin(t)|^2
  \non
  \\
  && \quad {}
  + \intQt |\nabla\dt\phin|^2
  + \frac 12 \iO |\nabla\phin(t)|^2
  + \intQt \betaeps'(\phin)|\dt\phin|^2
  + \intQt \phin \, \betaeps'(\phin) \dt\phin
  \non
  \\
  && = \frac 12 \iO |\nabla\mun(0)|^2
  + \frac {L'}2 \iO |\phin(0)|^2
  + \frac 12 \iO |\nabla\phin(0)|^2
  \non
  \\
  && \quad {}
  + \intQt u \, \dt\mun
  + L' \intQt u \, \dt\phin
  - L' \intQt \nabla\mun \cdot \nabla\dt\phin
  \non
  \\
  && \quad {}
  - \intQt \pi'(\phin) |\dt\phin|^2
  - \intQt \phin \, \pi'(\phin) \dt\phin \,.
  \label{perstimareg}
\Esist
All of the terms on the \lhs\ are nonnegative but the last one.
To treat it, we introduce $\gamma,\,\hat\gamma:\erre\to\erre$
by setting
\Beq
  \gamma(r) := r \, \pier{\betaeps}'(r)
  \aand
  \hat\gamma(r) := \int_0^r \gamma(s) \, ds
  \quad \hbox{for $r\in\erre$},
  \non
\Eeq
and notice that $\hat\gamma$ is nonnegative.
Hence, we have~that
\Beq
   \intQt \phin \, \betaeps'(\phin) \dt\phin
   = \iO \hat\gamma(\phin(t))
   - \iO \hat\gamma(\phin(0))
   \geq - \iO \hat\gamma(\phin(0)) \,.
   \non
\Eeq
On the other hand, by assuming that $n$ is large enough
as we did to obtain~\eqref{betaepsphinz}, and owing to~\eqref{stimabetaprimo}, 
we can estimate, \nju{at first $\norma{\betaeps'(\phin(0))}_\infty$ and then $\norma{\hat\gamma(\phin(0))}_\infty$}, 
uniformly with respect to $n$ and~$\eps$.
Let us come to the \rhs.
While the second and third \nju{terms} can easily be estimated by means of \eqref{normephinz},
the first one needs some work.
We write \eqref{secondan} at $t=0$ and test it by $-\Delta\mun(0)$.
We obtain that
\Beq
  \iO |\nabla\mun(0)|^2
  = \iO \bigl\{
    \nabla(-\Delta\phin(0))
    + \bigl( \betaeps'(\phin(0)) + \pi'(\phin(0) \bigr) \nabla\phin(0)
  \bigl\} \cdot \nabla\mun(0) \,.
  \non
\Eeq
By arguing as in Remark~\ref{Phinz}, we see that
$\norma{\nabla(-\Delta\phin(0))}\leq c\,\norma\phiz_{\Hx3}$.
Moreover, $\norma{\betaeps'(\phin(0))}_\infty$ has been already estimated,
and $\phin(0)$ is uniformly bounded in $V$ by~\eqref{normephinz}.
Therefore, by the Schwarz inequality \pier{it is clear that}
\Beq
  \norma{\nabla\mun(0)} \leq c \,.
  \label{nablamunz}
	\Eeq
Let us \nju{proceed} and treat the first integral involving~$u$. 
By integrating by parts, and then owing to the Poincar\'e 
\pier{and Young inequalities} and our assumption on~$u$, 
we have for every~$\delta>0$ that
\Bsist
  && \intQt u \, \dt\mun
  = \iO u(t) \mun(t) 
  - \iO u(0) \mun(0)
  - \intQt \dt u \, \mun
  \non
  \\
  \separa
  && = \iO u(t) (\mun(t)-\barmun(t))
  - \iO u(0) (\mun(0)-\barmun(0))
  - \intQt \dt u (\mun-\barmun)
  \non
  \\
  && \quad {}
  + \iO u(t) \barmun(t) 
  - \iO u(0) \barmun(0) 
  - \intQt \dt u \, \barmun
  \non
  \\
  \separa
  && \leq \delta \iO |\nabla\mun(t)|^2
  + \iO |\nabla\mun(0)|^2
  + \intQt |\nabla\mun|^2
  \non
  \\
  && \quad {}
  + M |\Omega| \, |\barmun(t)|
  + M |\Omega| \, |\barmun(0)|
  + |\Omega| \ioT |\barmun(s)|^2 \, ds
  + c_\delta \,.
  \label{uterm}
\Esist
Since \eqref{nablamunz} has been established, we just have to estimate $\barmun$ pointwise.
We observe that testing \pier{\eqref{secondan}} by $1/|\Omega|$ and using \eqref{secondastiman} yields that
\Beq
  |\barmun(t)|
  \leq |\Omega|^{-1} \norma{\betaeps(\phin(t))}_1 
  + c 
  \label{prebarmunt}
\Eeq
for every $t\in[0,T]$.
On the other hand, we can test \pier{\eqref{secondan}} by $\phin(t)-\barphin(t)$, 
and account for \eqref{fromMZ}, the Poincar\'e inequality, and~\eqref{secondastiman},
to have~that
\Bsist
  && \delta_0 \, \norma{\betaeps(\phin(t))}_1
  \leq \iO (\mun(t)-\barmun(t)) (\phin(t)-\barphin(t))
  - \iO \pi(\phin(t)) (\phin(t)-\barphin(t))
  + c
  \non
  \\
  && \leq c \, \norma{\nabla\mun(t)} 
  + c \,.
  \non
\Esist
By combining, we deduce that
\Beq
  |\barmun(t)|
  \leq c \, \norma{\nabla\mun(t)} 
  + c 
  \quad \hbox{for every $t\in[0,T]$}.
  \label{barmunt}
\Eeq
Hence, by also accounting \nju{for} \eqref{secondastiman} \pier{and \eqref{nablamunz}}, we have that
\Bsist
  && M |\Omega| \, |\barmun(t)|
  \leq \delta \iO |\nabla\mun(t)|^2 + c_\delta\,, 
  \non
  \\
  && |\barmun(0)|
  \leq c \iO |\nabla\mun(0)|^2 + c
  \leq c
  \aand
  \ioT |\barmun(s)|^2 \, ds
  \leq c \intQ |\nabla\mun|^2 + c
  \leq c \,,
  \non 
\Esist
and \eqref{uterm} yields
\Beq
  \intQt u \, \dt\mun
  \leq 2 \delta \iO |\nabla\mun(t)|^2
  + c_\delta \,.
  \non
\Eeq
Next, we have
\Beq
  L' \intQt u \, \dt\phin
  - L' \intQt \nabla\mun \cdot \nabla\dt\phin
  \leq \delta \intQt |\dt\phin|^2
  + \delta \intQt |\nabla\dt\phin|^2
  + c_\delta\,,
  \non
\Eeq
as well as
\Beq
  - \intQt \pi'(\phin) |\dt\phin|^2
  - \intQt \phin \, \pi'(\phin) \dt\phin
  \leq L \intQt |\dt\phin|^2
  + \delta \intQt |\dt\phin|^2
  + c_\delta \,.
  \non
\Eeq
By combining all these estimates and \eqref{perstimareg},
recalling that $L'=L+1$, and choosing $\delta$ small enough,
we conclude that
\Beq
  \norma{\nabla\mun}_{\L\infty H} + \pier{\norma{\dt\phin}_{\L2V}}
  \leq c \,.
  \label{stimareg}
\Eeq
By also accounting for \eqref{barmunt}, we have~that
\Beq
  \norma\mun_{\L\infty V} \leq c \,.
  \label{stimaregmun}
\Eeq

\step
Conclusion

At this point we can let $n$ \nju{tend} to infinity and conclude that
\Beq
  \norma\phieps_{\L\infty V}
  + \pier{\norma{\dt\phieps}_{\L2V}}
  + \norma\mueps_{\L\infty V}
  \leq c \,,
  \label{stimaregmueps}
\Eeq
where $(\phieps,\mueps)$ is the solution to the approximating problem \gianni{\Pbleps}
(see Section~\ref{ProofWPeps}).
Moreover, testing \eqref{secondaeps} by $\betaeps(\phieps(t))$ \aat,
and then using \pier{the elliptic regularity theory},
we easily infer that
\Beq
  \norma{\betaeps(\phieps(t))}
  \leq \norma{\mueps(t)-\pi(\phieps(t))}
  \leq c
  \aand
  \normaW{\phieps(t)}
  \leq c 
  \quad \aat,
  \non
\Eeq
whence 
\Beq
  \norma{\betaeps(\phieps)}_{\L\infty H}
  + \norma\phieps_{\L\infty W} 
  \leq c \,.
  \non
\Eeq
\gianni{Now, we conclude.
Indeed, by letting $\eps$ tend to zero along a proper subsequence,
we obtain that the solution $(\phi,\mu,\xi)$ given by \Convsoluzeps\
has the expected regularity and satisfies estimate~\eqref{regestimate}
with a constant $K_2$ as in the statement.}

\subsection{The separation property}

In this section, we prove our \gianni{results} regarding the separation property.
Our proof of Proposition \ref{Separation} is based on a simple \gianni{consideration} on an elliptic problem.
For a given $g\in H$, let us consider the problem of finding a pair $(z,\zeta)\in V\times H$ satisfying
\Beq
  \gianni{\iO z v + {}}
  \iO \nabla z \cdot \nabla v
  + \iO \zeta v
  = \iO g v
  \quad \hbox{for every $v\in V$}
  \aand
  \zeta \in \beta(z) \quad \aeO \,.
  \label{elliptic}
\Eeq
Since $\beta$ is maximal monotone, this problem has a unique solution,
and its solution is the \gianni{weak limit in $V\times H$} of $(\zeps,\betaeps(\zeps))$, 
where $\zeps\in V$ solves
\Beq
  \gianni{\iO \zeps v +{}}
  \iO \nabla\zeps \cdot \nabla v
  + \iO \betaeps(\zeps) v
  = \iO g v
  \quad \hbox{for every $v\in V$} \,.
  \label{elleps}
\Eeq
For the reader's convenience, we give a short detail \gianni{on our last \nju{claim}}.
\gianni{Since $\betaeps$ is monotone and $\betaeps(0)=0$,}
by testing \eqref{elleps} first by $\zeps$, and then by $\betaeps(\zeps)$, we immediately obtain~that
\Beq
  \normaV\zeps \leq \norma g
  \aand
  \norma{\betaeps(\zeps)} \leq \norma g \,.
  \non
\Eeq
Then, letting $\eps$ tend to zero,
we have weak convergence in $V$ and strong convergence in $H$ for $\zeps$,
and weak convergence in $H$ for $\betaeps(\zeps)$,
so that the limiting pair of $(\zeps,\betaeps(\zeps))$ solves \eqref{elliptic}.

\Blem
\label{Elliptic}
Assume that $g\in\Linfty$.
Then, 
\Beq
  \zeta\in\Linfty
  \aand
  \norma\zeta_\infty
  \leq \norma g_\infty \,.
  \label{tesilemma}
\Eeq
\Elem

\Bdim
For \gianni{$k\in\enne$}, we introduce \gianni{$\zk\in V\cap\Linfty$} by setting
\Beq
  \zk := \min\{ k , \max\{ \betaeps(\zeps) , -k \} \}
  \non
\Eeq
and take $v=|\zk|^{p-2}\zk$ in \eqref{elleps}, \gianni{where $p>2$ is arbitrary}.
Since \gianni{$\zeps\zk\geq0$, $\nabla\zeps\cdot\nabla\zk\geq0$ and $\betaeps(\zeps)\zk\geq|\zk|^2$}, 
by also applying the \Holder\ and Young inequalities, we obtain~that
\Bsist
   && \norma\zk_p^p
  \leq \gianni{\iO \zeps |\zk|^{p-2}\zk + {}}
  \iO \nabla\zeps \cdot \nabla(|\zk|^{p-2}\zk)
  + \iO \betaeps(\zeps) |\zk|^{p-2} \zk
  = \iO g |\zk|^{p-2}\zk
  \non
  \\
  && \leq \norma g_p \, \norma{|\zk|^{p-2}\zk}_{p'}
  \leq \frac 1p \, \norma g_p^p
  + \frac 1{p'} \, \norma{|\zk|^{p-2}\zk}_{p'}^{p'}
  = \frac 1p \, \norma g_p^p
  + \frac 1{p'} \norma\zk_p^p\,,
  \non
\Esist
whence $\norma\zk_p\leq\norma g_p$.
By letting $k$ tend to infinity, we deduce that $\,\betaeps(\zeps)$\, is bounded and that
$\norma{\betaeps(\zeps)}_\infty\leq\norma g_\infty$.
\gianni{Thus, \nju{\,$\{\betaeps(\zeps)\}$\,} has a weak star limit in $\Linfty$ as $\eps\searrow0$.
Since it converges to $\zeta$ weakly in $\Ldue$, we infer that this limit is~$\zeta$.
Therefore, \pier{due to the weak star} lower semicontinuity of the $L^\infty$ norm, 
we obtain~\eqref{tesilemma}.}
\Edim

\step
Proof of Proposition~\ref{Separation}

\gianni{We rearrange \eqref{seconda} and apply the lemma by choosing, \aat, 
$z=\phi(t)$ and $g=\phi(t)+\mu(t)-\pi(\phi(t)\betti{)}$,
so that $\zeta=\xi(t)$.
Hence
\Beq
  \xi(t) \in \Linfty
  \aand
  \norma{\xi(t)}_\infty \leq \norma{\phi(t)+\mu(t)-\pi(\phi(t)\betti{)}}_\infty
  \quad \aat \pier{,}
  \non
\Eeq
that is, \eqref{xibdd} \nju{is valid}.
The last \nju{claim} of the statement easily follows.
Indeed, if $D(\beta)$ is a bounded open interval, then by maximal monotonicity we have that
$\betaz$ tends to infinity at the end-points of~$D(\beta)$.
Hence, there exists a compact interval $[a,b]$ such that $|\betaz(r)|>\norma\xi_\infty$
for every $r\in D(\beta)\setminus[a,b]$.
Since $\xi\in\beta(\phi)$ \aeQ, we also have that $|\xi|\geq|\betaz(\phi)|$ \aeQ,
whence $\phi\in[a,b]$ \aeQ.}\QED

\bigskip

Let us come to the separation property in dimension~$d=2$, i.e., Proposition~\ref{Separation2d},
whose assumptions are understood to be in force.
In our proof, we \nju{employ} some ideas and methods \nju{from} \cite{GGM}, but our argument is rather different.
First, we need a new approximation $\betaeps$ of~$\beta$,
since the Yosida regularization \nju{does not seem to be} suitable for our purpose.

\step
New regularization

In the case of the logarithmic potential \eqref{logpot}, we can choose
\Beq
  \beta(r) = \ln \frac {1+r}{1-r}
  \quad \hbox{for $r\in D(\beta)=(-1,1)$},
  \non
\Eeq
and define $\Beta$ accordingly in order that $\Beta(0)=0$.
Since $\beta$ is odd, for $\eps\in(0,1)$ we define $\betaeps$ 
to be the odd function on~$\erre$ that satisfies
\Bsist
  && \betaeps(r) = \beta(r)
  \quad \hbox{for $0\leq r \pier{{}\leq{}} 1-\eps$},
  \non
  \\
  && \betaeps(r) = \beta(1-\eps) + \betaeps'(1-\eps) (r-(1-\eps))
  \quad \hbox{for $r>1-\eps$},
  \non
\Esist
so that $\betaeps$ is of class $C^1$ and \accorpa{monbetaeps}{disugbetaeps} are satisfied.
We observe that $\betaeps$ enjoys some \nju{additional} regularity:
indeed, $\betaeps'$ is piecewise \Lip\ continuous and globally continuous, thus \Lip\ continuous.
Moreover, if $\Betaeps$ is the primitive of $\betaeps$ that vanishes at the origin,
\eqref{disugBetaeps} is satisfied, too.
In the previous sections, we have used one \nju{further} property of the Yosida regularization.
Namely, we had to guarantee that \betti{$\xi\in\beta(\phi)$} by knowing that
$\phieps$ and $\betaeps(\phieps)$ converge to $\phi$ and $\xi$
strongly in $\LQ2$ and weakly in $\LQ2$, respectively.
This still holds for the new $\betaeps$, since
for every $v,w\in H$ with $w\in\beta(v)$ 
(i.e., $|v|<1$ and $w=\beta(v)$ \aeO) there exist $\{\veps\}$
such that $\veps$ and $\betaeps(\veps)$ strongly converge in $H$ 
to $v$ and~$w$, respectively.
Indeed, one can take 
$\veps := \min\{ 1-\eps , \max\{ \veps , -1+\eps \} \}$
and observe that the Lebesgue dominated convergence theorem can be applied
since it holds \aeO\ that 
$\veps\to v$, $\betaeps(\veps)\to\beta(v)$, $|\veps|\leq|v|$, and $|\betaeps(\veps)|\leq|\beta(v)|$.
All this ensures that the solution (still termed $(\phieps,\mueps)$, of course) 
to the \gianni{new approximating problem \Pbleps}\ 
converges to the unique solution to the original problem
in the sense of \accorpa{convphieps}{convxieps} as $\eps$ tends to zero.
Hence, new uniform estimates for the new $(\phieps,\mueps,\pier{\betaeps(\phieps)})$ provide corresponding estimates for~$(\phi,\mu,\xi)$.

\step
Two inequalities

We first prove that
\Beq
  \betaeps'(r) \leq 2 \, e^{\gianni{|\betaeps(r)|}}
  \quad \hbox{for every $r\in\erre$}.
  \label{betaepsok}
\Eeq
It suffices to consider positive values of~$r$.
If $r<1-\eps$\nju{, then} we have~that
\Beq
  \betaeps'(r) = \beta'(r)
  = \frac 2{(1-r)(1+r)}
  \leq 2 \, \frac {1+r}{1-r}
  = 2 \, e^{\beta(r)} 
  = 2 \, e^{\betaeps(r)}.
  \non
\Eeq
If $r\geq1-\eps$\nju{, then} we have that
\Beq
  \betaeps'(r)
  = \beta'(1-\eps)
  = \frac 2{\eps(2-\eps)}
  \leq \frac 2\eps
  = 2 \, e^{-\ln\eps}
  \leq 2 \, e^{\ln(2-\eps)-\ln\eps}
  = 2 \, e^{\beta(1-\eps)}
  \leq 2 \, e^{\betaeps(r)}.
  \non
\Eeq
Hence, \eqref{betaepsok} \nju{is valid}.

Now, we prove that, for every $p\geq1$,
there exist two positive constants $\kappa$ and $\kappa'$ such~that
\Beq
  rs \, e^{ps}
  \leq \frac 12 \, s^2 \, e^{ps} + e^{\kappa r} + \kappa'
  \quad \hbox{for every $r,s\geq0$}.
  \label{fromyoung}
\Eeq
\nju{To prove this claim, consider} the function $\psi:\erre\to(-\infty,+\infty]$
defined by
$\psi(s):=(1+s)\ln(1+s)-s$ if $s>-1$ and extended by $1$ and $+\infty$
at $s=-1$ and for $s<-1$, respectively.
Then, $\psi$ is convex and l.s.c.,  and a direct computation 
easily shows that its conjugate function is given by
$\psi^*(r)=e^r-r-1$ for $r\in\erre$.
Then, the Young inequality $r's'\leq\psi(s')+\psi^*(r')$ holds true for every $r',s'\in\erre$.
In particular, if $r,s\geq0$ and $\delta>0$, then we have~that
\Beq
  rs \, e^{ps}
  \leq \psi(\delta \, s e^{ps}) + \psi^*(r/\delta)
  \leq \psi(\delta \, s e^{ps}) + e^{r/\delta} \,.
  \non
\Eeq
On the other hand, we also have that
\Bsist
  && \psi(\delta s\,e^{ps})
  \leq (1 + \delta s\,e^{ps}) \ln (1 + \delta s\,e^{ps})
  \non
  \\
  \separa
  && \leq (1 + \delta s\,e^{ps}) \ln (e^{ps} + \delta s\,e^{ps})
  = (1 + \delta s\,e^{ps}) (ps + \ln (1 + \delta\,s))
  \non
  \\
  \separa
  && \leq (1 + \delta s\,e^{ps}) (ps + \delta s)
  = s (p+\delta)
  + \delta (p+\delta) s^2 \, e^{ps}
  \non
  \\
  \separa
  && \leq \delta \, s^2 + \frac {(p+\delta)^2}{4\delta}
  + \delta (p+\delta) s^2 e^{ps}
  \leq \delta (1+p+\delta) s^2 e^{ps}
  + \frac {(p+\delta)^2}{4\delta} \,.
  \non
\Esist
Then, \eqref{fromyoung} follows if we choose $\delta>0$ such that
$\delta (1+p+\delta)=1/2$, and set $\kappa:=1/\delta$ and
$\kappa':=(p+\delta)^2/(4\delta)$.

\step
The basic estimate

Here is the most important change with respect to~\cite{GGM},
since we come back to the discrete problem \gianni{\Pbln}\ 
instead of directly dealing with problem \Pbl.
We can take advantage of all of the estimates of Section~\ref{UNIFEST}.
First, we notice that $\betaeps+\pi$ is a $C^1$ function
and recall that $u\in\H1H$.
It follows that $\phin$ and $\mun$ are functions \pier{in $\H2H$}.
Then, we can differentiate both \eqref{priman} and \eqref{secondan} with respect to time
and test the resulting inequalities by $\dt\phin$ and $\Delta\dt\phin$, respectively.
If we sum up and integrate \pier{by parts and} over $(0,t)$, then a cancellation occurs, and we obtain~that
\Bsist
  && \frac 12 \iO |\dt\phin(t)|^2
  + \intQt |\dt\phin|^2
  + \intQt |\Delta\dt\phin|^2
  \non
  \\
  && = \frac 12 \iO |\dt\phin(0)|^2
  + \intQt \dt u \, \dt\phin
  + \intQt (\betaeps'+\pi')(\phin) \, \dt\phin \, \Delta\dt\phin \,.
  \label{forbasic}
\Esist
We estimate the first term on the \rhs\ later on.
The next term can be treated in an obvious way by the Young inequality,
while the last one needs some work.
We have~that
\Bsist
  && \intQt (\betaeps'+\pi')(\phin) \, \dt\phin \, \Delta\dt\phin
  \leq \iot \norma{(\betaeps'+\pi')(\phin(s))}_3 \, \norma{\dt\phin(s)}_6 \, \norma{\Delta\dt\phin(s)} \, ds
  \non
  \\
  && \leq \frac 14 \intQt |\Delta\dt\phin|^2
  + \iot \norma{(\betaeps'+\pi')(\phin(s))}_3^2 \, \norma{\dt\phin(s)}_6^2 \, ds\,,
  \non
\Esist
and we have to estimate the last integral.
We have \aet~that
\Bsist
  && \norma{\dt\phin}_6^2
  \leq c \, \normaV{\dt\phin}^2
  \leq c \, \Bigl( |\dt\barphin|^2 + \iO |\nabla\dt\phin|^2 \Bigr)
  \non
  \\
  && \leq c + c \Bigl| \iO \dt\phin (-\Delta\dt\phin) \Bigr|
  \leq c + c \, \norma{\dt\phin} \, \norma{\Delta\dt\phin}\,,
  \non
\Esist
and we deduce that
\Bsist
  && \iot \norma{(\betaeps'+\pi')(\phin(s))}_3^2 \, \norma{\dt\phin(s)}_6^2 \, ds
  \non
  \\
  && \leq c \iot \norma{(\betaeps'+\pi')(\phin(s))}_3^2 \,
      \bigl( c + c \, \norma{\dt\phin(s)} \, \norma{\Delta\dt\phin(s)} \bigr) \, ds
  \non
  \\
  \separa
  && \leq c \iot \norma{(\betaeps'+\pi')(\phin(s))}_3^2 \, ds
  \non
  \\
  && \quad {}
  + \frac 14 \intQt |\Delta\dt\phin|^2
  + c \iot \norma{(\betaeps'+\pi')(\phin(s))}_3^4 \, \norma{\dt\phin(s)}^2 \, ds
  \non
  \\
  \separa
  && \leq \betti{c\iot \left(\norma{\betaeps'(\phin(s))}_3^2+1\right)} \, ds
  \non
  \\
  && \quad {}+ \frac 14 \intQt |\Delta\dt\phin|^2
  + c \iot \bigl( 1 + \norma{\betaeps'(\phin(s))}_3^4 \bigr) \norma{\dt\phin(s)}^2 \, ds \,.
  \non
\Esist
By combining these estimates with \eqref{forbasic} and applying the Gronwall lemma,
we conclude~that
\Bsist
  && \norma{\dt\phin}_{\L\infty H}^2
  + \norma{\Delta\dt\phin}_{\L2H}^2
  \non
  \\[2mm]
  && \leq c \, \bigl(
    \norma{\dt\phin(0)}^2
    + \norma{\dt u}_{\L2H}^2
    + \norma{\betaeps'(\phin)}_{\L2{\Lx3}}^2 + 1
  \bigr) \times {}
  \non
  \\
  && \quad {}
  \times \exp \Bigl( c \ioT \bigl( 1 + \norma{\betaeps'(\phin(s)}_3^4 \bigl) ds \Bigr) \,,
  \non
\Esist
whence also
\Bsist
  && \norma{\dt\phin}_{\L\infty H}^2
  + \norma{\Delta\dt\phin}_{\L2H}^2
  \non
  \\[2mm]
  && \leq c \, \bigl(
    \norma{\dt\phin(0)}^2
    + \norma{\betaeps'(\phin)}_{\pier{\L2{\Lx3}}}^2 + 1
  \bigr) 
  \, e^{c \, \norma{\betaeps'(\phin)}_{\pier{\L4{\Lx3}}}^4} \,.
  \label{basicn}
\Esist
In order to let $n$ tend to infinity, we have to estimate the $H$ norm of~$\dt\phin(0)$.
To this end, we set $\eps_0:=(1-\norma\phiz_\infty)/2$, assume that $\eps\leq\eps_0$,
and recall that $\phin(0)$ converges uniformly to~$\phiz$ 
(see also Remark~\ref{Phinz}). 
Hence, we can assume $n$ large enough
in order that $\norma{\phin(0)}\leq1-\eps_0$.
On the other hand, $(\betaeps+\pi)(r)=f(r)$, the whole logarithmic potential,
for $|r|\leq1-\eps_0$, and $f$ is smooth in~$(-1,1)$.
Therefore, $(\betaeps+\pi)(\phin(0))=f(\phin(0))$ is as smooth as~$\phin(0)$
(i.e., as the eigenfunctions~$e_j$),
and, for $0\leq s\leq4$, the norm of $f(\phin(0))$ in $\Hx s$ can be uniformly estimated
by the corresponding norm of~$\phin(0)$, thus by the norm of~$\phiz$, 
due to our further assumption~\eqref{newhpphiz}
(see Remark~\ref{Phinz} once more). 
At this point, we start estimating.
We test \eqref{priman} and \eqref{secondan}, written with $t=0$, 
by $\dt\phin(0)$ and $\Delta\dt\phin(0)$, respectively, and sum up.
Since the terms involving $\mun(0)$ cancel each other, \pier{with the help of some integrations by parts} we obtain~that
\Bsist
  && \norma{\dt\phin(0)}^2
  = \iO |\dt\phin(0)|^2
  \non
  \\
  && = \iO \bigl( u(0) - \phin(0) \bigr) \, \dt\phin(0)
  + \iO \nabla\phin(0) \cdot \nabla\Delta\dt\phin(0)
  + \iO f(\phin(0)\pier) \Delta\dt\phin(0)
  \non
  \\
  && = \iO \bigl\{
      u(0) - \phin(0) 
      - \Delta^2 \phin(0)
      + \Delta (f(\phin(0))
  \bigr\} \dt\phin(0)  
  \non
  \\[1mm]
  && \leq c \, \norma{\dt\phin(0)} \,.
  \non
\Esist
At this point, we are ready to let $n$ tend to infinity in~\eqref{basicn}.
Indeed, by recalling \eqref{convphin} and applying, e.g, \cite[Sect.~8, Cor.~4]{Simon}), 
we have that $\phin$ converges to $\phieps$ strongly in $\C0{\Lx3}$.
Since $\betaeps'$ is \Lip\ continuous, this implies that 
$\betaeps'(\phin)$ converges to $\betaeps'(\phieps)$ in the same topology.
As a consequence, the norm of $\betaeps'(\phin)$ in $\L p{\Lx3}$
converges to \nju{the} corresponding norm of $\betaeps'(\phieps)$ for $p=2$ and $p=4$,
so that, taking the limit in \eqref{basicn}, we obtain an inequality.
Then, by estimating the above norms with the norm in $\L\infty{\Lx3}$,
we conclude that
\Bsist
  && \norma{\dt\phieps}_{\L\infty H}^2
  + \norma{\Delta\dt\phieps}_{\L2H}^2
  \non
  \\[2mm]
  && \leq c \, \bigl(
    1
    + \norma{\betaeps'(\phieps)}_{\L\infty{\Lx3}}^2
  \bigr) 
   \, e^{c \, \norma{\betaeps'(\phieps)}_{\L\infty{\Lx3}}^4} \,.
  \label{basic}
\Esist

\step
Conclusion of the proof \gianni{of Proposition~\ref{Separation2d}}

Our aim is to derive a bound in $\L\infty H$ for $\dt\phi$.
To this end, we estimate the \rhs\ of~\eqref{basic}.
Here, we follow \cite{GGM} rather closely.
We set $\phiepsk := \min\{ k , \max\{ \phieps , -k \} \}$
and $\Psieps(r):=\betaeps(r)e^{3|\pier{\betaeps(r)}|}$ for $r\in\erre$
and observe that $v=\Psieps(\phiepsk(t))$ is admissible in \eqref{secondaeps} \aat.
Therefore, \aet\ we have~that
\Beq
  \iO \Psieps'(\phiepsk) \nabla\phieps \cdot \nabla\phiepsk
  + \iO \betaeps(\phieps) \Psieps(\phiepsk)
  = \iO \geps \Psieps(\phiepsk) \,,
  \non
\Eeq
where, for brevity, we have set
$\,\geps=\mueps-\pi(\phieps)$.
Since $\Psieps'$ is nonnegative, $\nabla\phieps\cdot\nabla\phiepsk=|\nabla\phiepsk|^2$ a.e.,
and $\betaeps(\phieps)\betaeps(\phiepsk)\geq|\betaeps(\phiepsk)|^2$ a.e.,
we deduce~that
\Beq
  \iO |\betaeps(\phiepsk)|^2 \, e^{3|\betaeps(\phiepsk)|}
  \leq \iO \geps \Psieps(\phiepsk) \,.
  \non 
\Eeq
We estimate the \rhs\ \pier{using the inequality} \eqref{fromyoung} with $p=3$, \betti{$r=|\geps|$, and $s=|\betaeps(\phiepsk)|$} and have that
\Beq
  \iO \geps \Psieps(\phiepsk)
  \leq \iO |\geps| \, |\betaeps(\phiepsk)| e^{3|\betaeps(\phiepsk)|}
  \leq \frac 12  \iO |\betaeps(\phiepsk)|^2 \, e^{3|\betaeps(\phiepsk)|}
  + \iO \bigl( e^{\kappa|\geps|} + \kappa' \bigr) .
  \non
\Eeq
On the other hand, by observing that $(1-r^2)e^r\leq e$ for every $r\in\erre$,
we also have~that
\Beq
  \iO e^{3|\betaeps(\phiepsk)|}
  \leq 9 \iO |\betaeps(\phiepsk)|^2 e^{3|\betaeps(\phiepsk)|}
  + c \,.
  \non
\Eeq
By combining all this with \eqref{betaepsok}, we infer that
\Beq
  \iO |\betaeps'(\phiepsk)|^3
  \leq c \iO e^{3|\betaeps(\phiepsk)|}
  \leq c \iO e^{\kappa|\geps|}
  + c \,.
  \non
\Eeq
At this point, we recall the Trudinger inequality \pier{(see, e.g., \cite{Moser})}

\Beq
  \iO e^{|v|}
  \leq \CO \, e^{\CO\,\normaV v^2}
  \quad \hbox{for every $v\in V$}\,,
  \non
\Eeq
\pier{which holds since we are supposing that $d=2$, and} where the constant $\CO$ only depends on~$\Omega$.
Hence, we conclude that
\Beq
  \iO |\betaeps'(\phiepsk)|^3
  \leq c \iO e^{\kappa|\geps|}
  + c
  \leq c \, e^{\kappa^2 \pier{\CO} \normaV{\mueps-\pi(\phieps)}^2}
  \quad \aet \,.
  \non
\Eeq
By accounting for \eqref{fromunifest} and~\eqref{stimaregmueps},
\pier{in addition we infer} that
\Beq
  \norma{\betaeps'(\phiepsk)}_{\L\infty{\Lx3}} 
  \leq c \,,
  \quad \hbox{whence also} \quad
  \norma{\betaeps'(\phieps)}_{\L\infty{\Lx3}} 
  \leq c \,.
  \non
\Eeq
At this point, we can take the limit in \eqref{basic} and deduce that
$\dt\phi$ and $\Delta\dt\phi$ belong to $\L\infty H$ and $\L2H$, respectively.
Therefore, $\dt\phi$ belongs to $\L2W$.
Moreover, elliptic regularity in \eqref{prima} implies that
$\mu$ belongs to $\L\infty W\subset\LQ\infty$.
Since $\pi(\phi)$ is bounded, too, by applying Proposition~\ref{Separation},
we obtain that $\xi$ is bounded.
This concludes the proof of~\eqref{reg2d}.
Furthermore, it is clear that \eqref{separation2d} holds true as well,
with a constant $K_3$ as in the statement.
Thus, the proof of Proposition~\ref{Separation2d} is complete.


\section{The control problem}
\label{CONTROL}
\setcounter{equation}{0}

In this section, we investigate the optimal control problem \eqref{controlPbl}.
We remark that we will use some notations already utilized in the previous sections with a different meaning
(e.g.,~\gianni{$\xi$} in the linearized system introduced later~on).
However, no confusion can arise.

\subsection{Proof of Theorem~\ref{Optimum}}
\juerg{We assume first that $\alpha_3=0$ in which case the cost functional is independent of the
solution variables $\mu,\xi$. 
By Theorem~\gianni{\ref{Wellposedness}}, the mapping $\Uad\ni u\mapsto \phi\,$ is well
defined, and thus also the cost functional. 
Now we pick any minimizing sequence $\{((\phi_n,\mu_n,\xi_n),u_n)\}$
for the optimal control problem, that is, in particular, $(\phi_n,\mu_n,\xi_n)$,
where $\xin\in\beta(\phin)$, satisfies \Pbl\ with 
\rhs\ $\,u=u_n$ for all $n\in\enne$. 
Owing to \pier{\eqref{defUad} and}  the global estimate \eqref{stimaK1}, we may without loss of generality
assume that there are limits $\phi,\mu,\xi,u$ such that}
\begin{align}
\label{exi1}
u_n\to u&\quad\mbox{weakly star in $\,\pier{{}\H1H\cap{}}\LQ\infty$},\\
\label{exi2}
\phin\to\phi&\quad\mbox{weakly star in }\,H^1(0,T;V^*)\cap L^\infty(0,T;V)\cap L^2(0,T;W)\nonumber\\
&\quad\mbox{and strongly in }\,C^0([0,T];H),\\
\label{exi3}
\mun\to\mu&\quad\mbox{weakly in }\,L^2(0,T;V),\\
\label{exi4}
\xin\to\xi&\quad\mbox{weakly in }\,L^2(0,T;H),
\end{align}
\juerg{where the strong convergence in \eqref{exi2} follows from, e.g., \cite[Sect.~8,~Cor.~4]{Simon}.
Obviously, we have that $u\in\Uad$. 
Now, \gianni{by also recalling Remark~\ref{IntPbl}}, 
we can pass to the limit as $n\to\infty$ in the state system \Pbl, written with \rhs\ $u_n$, to see
that $(\phi,\mu,\xi)$ solves \Pbl\ with \rhs~$\,u$. 
In view of \eqref{exi4} and the strong
convergence stated in \eqref{exi2}, and since \pier{the extension of $\beta$ to $L^2(Q)$} is maximal monotone, a standard argument 
then yields that $\xi\in\beta(\phi)$. The pair $((\phi,\mu,\xi),u)$ is thus admissible for \eqref{controlPbl},
and the weak sequential lower semicontinuity properties of the cost functional yield that it is a
minimizer.}   

\juerg{
Suppose now that $\beta$ is single-valued, in which case the cost functional is well defined also
if $\alpha_3>0$. 
Then the above line of argumentation can obviously be repeated, only that $\xi_n=\beta(\phin)$
and $\xi=\beta(\phi)$ in this case. 
The assertion is thus proved.\QED  }

\subsection{Fr\'echet differentiability of the control-to-state operator}

In this section, we show the Fr\'echet differentiability of the control-to-state $\S$
introduced in \eqref{defS}. 
To this end, we recall the definition of $\UR$ 
given in \eqref{defUR} and the assumptions~\Start\ which are assumed to hold throughout 
the remainder of this section. 
Observe that they are fulfilled in the case $d=3$ for 
everywhere defined smooth potentials and in the case $d=2$ for the logarithmic 
potential \eqref{logpot} under the assumptions of Proposition~\gianni{\ref{Separation2d}}.     

Now let some $\ustar\in\UR$ be fixed and $(\phistar,\mustar)=\calS(\ustar)$ be 
the associated state. 
We then consider the corresponding linearized system,
that is, given $h\in {\cal X}$, we seek  a pair $(\xi,\eta)$ satisfying
\Bsist
  && \xi \in \H1\Vp \cap \L\infty V \cap \L2W
  \aand
  \eta \in \L2V\,,
  \qquad
  \label{regxieta}
  \\
  && \< \dt\xi(t) , v > 
  + \iO \xi(t) v 
  + \iO \nabla\eta(t) \cdot \nabla v
  = \iO h(t) v
  \non
  \\
  && \quad \hbox{\aat\ and every $v\in V$}\,,
  \label{primaL}
  \\
  && \iO \eta(t) v 
  = \iO \nabla\xi(t) \cdot \nabla v
  + \iO f''(\phistar(t)) \xi(t) v
  \non
  \\
  && \quad \hbox{\aat\ and every $v\in V$}\,,
  \label{secondaL}
  \\
  && \xi(0) = 0 \,.
  \label{cauchyL}
\Esist
\Accorpa\PblL primaL cauchyL

We have the following result.

\Blem
\label{WellposednessL}
Suppose that the conditions \Start\ are fulfilled, and let $\ustar\in\UR$ be given with
associated state $(\phistar,\mustar)=\S(u^*)$. 
Then the system \PblL\ has for every $h\in{\cal X}$
a unique solution $(\xi,\eta)$ with the regularity \eqref{regxieta}. 
Moreover, the linear mapping
$h\mapsto(\xi,\eta)$ is continuous from ${\cal X}$ into $\calY$.
\Elem

\Bdim 
The existence proof is similar to that for the state system, 
and we just provide a sketch, leaving the details to the reader. 
Indeed, the system \PblL\ has the same structure for $(\xi,\eta)$ as the state system \Pbl\ 
for $(\phi,\mu)$, only that we have zero initial conditions here and the term $f''(\phistar)\xi$ in place
of $\beta(\phi)+\pi(\phi)$ (recall that $\beta$ is single-valued). 
We thus may argue as in Section~\gianni{\ref{EXISTENCE}}:
one approximates the system \PblL\ by a Faedo--Galerkin scheme using
the eigenfunctions \eqref{eigen} and the subspaces \eqref{defVn}. 
Since the term $f''(\phistar)\xi$,
where $f''(\phistar)\in L^\infty(Q)$, is much easier to handle than the nonlinearities $\beta$ and~$\pi$,
similar, but in comparison with the state system considerably simpler, a~priori estimates can be performed 
on the Faedo--Galerkin system, yielding the bound (compare \eqref{stimaK1})
$$
\|\xi_n\|_{H^1(0,T;V^*)\cap L^\infty(0,T;V)\cap L^2(0,T;W)}\,+\,\|\eta_n\|_{L^2(0,T;V)}\,\le\,c 
\quad\forall\,n\in\enne.
$$
Using compactness arguments, we can pass to the limit as $n\to\infty$ in the Faedo--Galerkin system
on a subsequence in order to establish the existence result.

We only prove that the mapping $h\mapsto (\xi,\eta)$ has the asserted continuity property. 
To this end, let $h\in {\cal X}$ be fixed and $(\xi,\eta)$ be an associated solution to \PblL\ 
with the regularity \eqref{regxieta}. 
We then insert $\betti{v=}\xi$ in \eqref{primaL} and $\betti{v=}\eta$ in \eqref{secondaL}, 
add the resulting identities, and integrate over $(0,t)$ for arbitrary $t\in (0,T]$. 
After a cancellation of two terms, we then obtain that 
\begin{align}
&\frac 12\,\|\xi(t)\|^2+\intQt |\xi|^2+\intQt |\eta|^2\,=\,\intQt h\,\xi
-\intQt f''(\phistar)\,\xi\,\eta\non\\
&\le\,\frac 12\intQt |\eta|^2\,+\,c\intQt (|\xi|^2+|h|^2)\,.\non
\end{align}
Now observe that, by virtue of continuous embedding, we obviously have that $\xi\in C^0([0,T];H)$. 
Therefore, Gronwall's lemma yields that
\begin{equation}
\label{conny1}
\|\xi\|_{C^0([0,t];H)}+\|\eta\|_{L^2(0,t;H)}\,\le\,c\,\|h\|_{L^2(0,t;H)}
\quad\mbox{for all $\,t\in [0,T]$}.
\end{equation}
Finally, applying elliptic regularity theory to \eqref{secondaL}, we 
may also conclude that
\begin{equation}
\label{conny3}
\|\xi\|_{L^2(0,t;W)}\,\le\,c\,\|h\|_{L^2(0,t;H)}
\quad\mbox{for all $\,t\in [0,T]$}.
\end{equation}
With this, the asserted continuity property is shown. 
Finally, it is easily seen that the inequality \eqref{conny1} 
also implies the uniqueness of the solution: 
indeed, \gianni{by using linearity, if $(\xi,\eta)$ solves \PblL\ with $h=0$, then
\eqref{conny1} implies that $\xi=\eta=0$}. 
\Edim

\Brem 
The existence proof sketched above yields that a unique solution to \PblL\ 
with the regularity \eqref{regxieta} also exists 
if we only have $h\in L^2(Q)$ and that also in this situation the continuity property shown above is valid.
\Erem

We now show the following \Frechet\ differentiability result.

\Bthm
\label{Frechet}
Suppose that the conditions \Start\ are fulfilled, and let $\ustar\in\UR$ be given with
associated state $(\phistar,\mustar)=\S (u^*)$.  
Then the operator $\calS$ is \Frechet\ differentiable at $\ustar$ as a mapping from $\calX$ into $\calY$, 
and the \Frechet\ derivative $D\calS(\ustar)\in {\cal L}(\calX,\cal Y)$ acts as follows:
for every  $h\in\calX$,  the value $D\calS(\ustar)(h)$ is given by
the solution $(\xi,\eta)$ to the linearized system \PblL.
\Ethm

\Bdim
Since $\ustar\in \UR$, it follows that there is some $\Lambda>0$ such that
$\ustar+h\in\UR$ whenever $\|h\|_{\cal X}\,\le \Lambda$. 
In the following, we only consider such perturbations $h\in\calX$. 
In the remainder of this proof, we denote by the small-case symbol $\,c\,$ constants that
may depend on the data of the system but not on the choice of $h\in\calX$ with 
$\,\|h\|_{\calX}\le\Lambda$.
We then define 
\begin{align}
&(\phi^h,\mu^h):={\cal S}(\ustar+h),\quad y^h:=\phi^h-\phistar-\xi^h,\quad z^h:=\mu^h-\mustar-\eta^h,
\end{align}
where $(\xi^h,\eta^h)$ denotes the unique solution to the linearized system \PblL. 
We then have
\begin{equation}\label{regyhzh}
y^h\in H^1(0,T;V^*)\cap L^\infty(0,T;V)\cap L^2(0,T;W),\quad z^h\in L^2(0,T;V).
\end{equation}

Now, by virtue of Lemma~\gianni{\ref{WellposednessL}}, 
the linear mapping $h\mapsto (\xi^h,\eta^h)$ is continuous between $\calX$ and~$\calY$. 
According to the definition of \Frechet\ differentiability, it therefore suffices
to construct a function $G:(0,\Lambda)\to (0,\infty)$ that satisfies 
$\lim_{\lambda\searrow0}G(\lambda)/(\lambda^2)=0$ and
\begin{equation}
\label{glambda}
\|(y^h,z^h)\|_{\calY}^2\,\le\,G(\|h\|_\infty)\,.
\end{equation}
We are going to show that we may choose $G(\lambda)=\hat c\,\lambda^4$ with some sufficiently
large $\hat c>0$.

Next, it is easily seen, using Taylor's theorem with integral remainder, that $(y^h,z^h)$ is a solution to the system
\begin{align}
&\langle \dt y^h(t),v\rangle + (y^h(t),v)+(\nabla z^h(t),\nabla v)=0, \quad\mbox{\aat\ and all
$v\in V$},
\label{yh1}
\\[1mm]
&(z^h(t),v)=(\nabla y^h(t),\nabla v)+(f''(\phistar(t))\,y^h(t),v)
+(R^h(t)(\phi^h(t)-\phistar(t))^2,v),\non\\
&\quad\mbox{\aat\ and all $v\in V$},
\label{yh2}
\\[1mm]
&\gianni{y^h(0)=0} \quad\mbox{a.e. in }\,\Omega,
\label{yh3} 
\end{align}
with the remainder
$$
  \gianni{R^h(t):= \int_0^1 (1-s)\,f'''(\phistar(t)+\,s\,(\phi^h(t)-\phistar(t))) \,ds\,.} 
$$
Observe that \gianni{$\phistar+s(\phi^h-\phistar)\in [a,b]$ almost everywhere in $Q$ for all $s\in [0,1]$}, and~thus
\begin{equation}
\|R^h\|_\infty \,\le\,c\quad\mbox{for all $h\in\calX$ satisfying $\|h\|_{\calX}\le\Lambda$.}
\label{Rhbound}
\end{equation}
Also, we conclude from \eqref{contdep3} in Theorem~\gianni{\ref{Lipcontdep}} that the estimate
\begin{equation}
\|\phi^h-\phistar\|_{C^0([0,t];H)\cap L^2(0,t;W)}\,+\,\|\mu^h-\mustar\|_{L^2(0,t;H)}
\,\le\,C_2\,\|h\|_{\pier{L^2(0,T;H)}}
\label{Ringo}
\end{equation}
is valid for all $t\in (0,T]$ and all admissible perturbations~$h$.
Furthermore, taking $v=1$ in \eqref{yh1} and accounting for \eqref{yh3}, we immediately see~that 
\begin{equation}\label{meanvyh}
\nju{\overline {y^h}}(t) =0\quad\mbox{for all $t\in[0,T]$}.
\end{equation} 
Therefore, we may insert $v={\cal N}(y^h)(t)$ in \eqref{yh1} and $v=-y^h(t)$ in \eqref{yh2},
add the resulting identities, and integrate over time. 
Noting an obvious cancellation, we then obtain that
\begin{align}
\label{John}
&\frac 12 \,\|y^h(t)\|_*^2\,+\int_0^t \|y^h(s)\|_*^2\,ds 
\,+\,\intQt |\nabla y^h|^2\non\\
&=\,-\intQt f''(\phistar)\,|y^h|^2\,-\,\intQt R^h\,(\phi^h-\phistar)^2\,y^h
\,=: \,I_1+I_2,
\end{align}
with obvious meaning. 
Now, owing to \eqref{globalbdd} and by virtue of the compactness inequality \eqref{compact}, we have that
\begin{align}
I_1\,&\le\,c\int_0^t\|y^h(s)\|^2\,ds\,\le\, \frac 14 \intQt |\nabla y^h|^2\,+\,c\int_0^t
\|y^h(s)\|_*^2\,ds\,.
\label{Paul}
\end{align}
Moreover, by also using Young's and H\"older's inequalities, 
\eqref{Ringo}, and \eqref{Rhbound}, we see that
\begin{align}
I_2\,&\le\,c\int_0^t \|(\phi^h-\phistar)(s)\|_4\,\|(\phi^h-\phistar)(s)\|_2\,
\|y^h(s)\|_4\,ds
\nonumber\\
&\le\,c\,\|\phi^h-\phistar\|_{C^0([0,t];H)}\int_0^t \|(\phi^h-\phistar)(s)\|_V\,
\|y^h(s)\|_V\,ds\nonumber
\\
&\le\,c\,\|h\|_{L^2(0,t;H)}\,\|\phi^h-\phistar\|_{L^2(0,t;V)}\,\gianni{\|y^h\|_{L^2(0,t;V)}}\non\\
&\le\,\frac 14\intQt|\nabla y^h|^2\,+c\int_0^t\|y^h(s)\|_*^2\,ds 
\,+\,c \|h\|^4_{\pier{L^2(0,T;H)}}\,.
\label{George}
\end{align}
\pier{Combining \eqref{John}--\eqref{George} with Gronwall's lemma, 
we thus have shown that}
\begin{equation}
\label{esti1}
\nju{\|y^h\|}^2_{L^\infty(0,t;V^*)\cap L^2(0,t;V)}\,\le\,c\,\|h\|^4_{\pier{L^2(0,T;H)}}\quad\mbox{for all $\,t\in[0,T]$}.
\end{equation}

The next step is to insert $v=y^h(t)$ in \eqref{yh1} and $v=z^h(t)$ in \eqref{yh2}, add the resulting equations,
and integrate over time. 
This yields the identity
\begin{align}
&\frac 12\,\|y^h(t)\|^2\,+\intQt |y^h|^2\,+\intQt |z^h|^2\non\\
&=\,\intQt f''(\phistar)\,y^h\,z^h\,+\intQt R^h\,(\phi^h-\phistar)^2\,z^h\,:=\,I_3+I_4\,.
\label{Mick}
\end{align}
Clearly, by \eqref{globalbdd} and Young's inequality \pier{we have that} 
\begin{equation}
I_3\,\le\,\frac 14 \intQt|z^h|^2\,+\,c\intQt |y^h|^2\,.
\label{Keith}
\end{equation}
Moreover, also using \eqref{Rhbound}, \eqref{Ringo}, and H\"older's inequality, we find that
\begin{align}
I_{\betti{4}}\,&\le\,c\int_0^t \|(\phi^h-\phistar)(s)\|_\infty\,\|(\phi^h-\phistar)(s)\|_2\,\|z^h(s)\|_2\,ds\non\\
&\le\,\frac 14\intQt |z^h|^2\,+\,c\,\|\phi^h-\phistar\|^2_{C^0([0,t];H)}
\int_0^t\|(\phi^h-\phistar)(s)\|_W^2\,ds \non\\
&\le\,\frac 14 \intQt |z^h|^2\, +\,c\,\|h\|^4_{\pier{L^2(0,T;H)}}\,.
\label{Charlie}
\end{align}
Now recall that $y^h\in C^0([0,T];H)$. 
Thus, combining \eqref{Mick}--\eqref{Charlie} with Gronwall's lemma,
and invoking \eqref{esti1}, we have finally shown that
$$
\|y^h\|_{C^0([0,T];H)\cap L^2(0,T;V)}^2\,+\,\gianni{\|z^h\|_{L^2(0,T;H)}^2}\, \le\,c_0\,\|h\|^4_{L^2(0,T;H)}\,,
$$
with a sufficiently large constant $c_0>0$. 
Since $\|h\|_{L^2(0,T;H)}\,\le\,T^{1/2}\,|\Omega|^{1/2}\,
\|h\|_\infty$, the condition \eqref{glambda} is fulfilled with 
\pier{$G(\lambda)= \hat c \lambda^4$ and}
$\,\hat c=T^2\,|\Omega|^2\,c_0$.
The assertion is thus proved.
\Edim

\subsection{Necessary conditions for optimality}

In this section, we establish first-order necessary optimality conditions 
for the optimal control problem. 
Using Theorem~\ref{Frechet} and the chain rule, a standard argument
yields the following result. 

\Bprop
\label{Badoptimality}
Suppose that the conditions \Start\ are fulfilled, and let $\ustar\in\Uad$ be an optimal
control with
associated state $(\phistar,\mustar)={\cal S} (u^*)$. Then it holds the variational inequality
\begin{align}
\label{vug1}
&\alpha_1\int_Q (\phistar-\phi_Q)\xi \,+\,\alpha_2\iO(\phistar(T)-\phi_\Omega)\xi(T)\,+
\,\alpha_3\int_Q(\mustar-\mu_Q)\eta\non\\
&+\,\alpha_4\int_Q u^*(u-u^*)\,\ge\,0\quad\forall\,u\in\Uad\,,
\end{align}
where $(\xi,\eta)$ is the unique solution to the linearized system \PblL\ with $h=u-\ustar$.
\Eprop

As usual, we now eliminate the variables $(\xi,\eta)$ from \eqref{vug1} by means of
the solution $(p,q)$ to the associated adjoint system, which is given by \Adjoint. 
We have the following result.

\Blem
Suppose that the conditions \Start\ are fulfilled, let $\ustar\in\Uad$ be an optimal
control with associated state $(\phistar,\mustar)={\cal S} (u^*)$, and assume that the regularity
conditions \eqref{ende}--\eqref{finale} are satisfied.
Then the adjoint system \Adjoint\ has a unique solution $(p,q)$ such that
\Beq
\label{regadjoint}
p\in H^1(0,T;V^*)\cap L^\infty(0,T;V)\cap \pier{{}L^2(0,T;W){}}, \quad q\in L^2(0,T;V).
\Eeq
\Elem

\Bdim
Again, we employ a Faedo--Galerkin approximation with the eigenfunctions defined in 
\eqref{eigen} and the $n$-dimensional subspaces $V_n$ from \eqref{defVn}. We then look
for functions 
$$\pier{p_n(x,t)=\sum_{i=1}^n p_n^i(t)e_i(x),\quad q_n(x,t)=\sum_{i=1}^n 
q_n^i(t)e_i(x),}$$
that solve the final value problem
\begin{align}
&(-\dt p_n(t),v)+(p_n(t),v)+(\nabla q_n(t),\nabla v)\,=\,-(f''(\phistar(t))q_n(t),v)
+(g_1(t),v)\non\\
&\qquad\mbox{\aat\ and every $v\in V_n$\,,}
\label{adn1}
\\
&(q_n(t),v)\,=\,(\nabla p_n(t),\nabla v)-(g_3(t),v)\non\\
&\qquad\mbox{\aat\ and every } v\in V_n\,,
\label{adn2}
\\
&(p_n(T),v)=(g_2,v)\quad\mbox{for every } v\in V_n\,,
\label{adn3}
\end{align}
where, owing to the regularity of $(\phistar,\mustar)$ and \eqref{ende}, 
\begin{eqnarray}
&&g_1=\alpha_1(\phistar-\phi_Q)\in L^2(0,T;V), \quad
g_2=\alpha_2(\phistar(T)-\phi_\Omega) \in V,\non\\
&&g_3=\alpha_3(\mustar-\mu_Q)\in L^2(0,T;V).
\label{paul}
\end{eqnarray}
Now observe that, by \pier{taking
$v= e_i$ in \eqref{adn2} and recalling \eqref{eigen}, we have that}
\begin{align}
\pier{q_n^i (t) =\lambda_i p_n^i (t) - (g_3(t),e_i)}
\quad\mbox{for $1\le i\le n$}.\non
\end{align} 
Therefore, if we insert $v=e_i$, $1\le i\le n$, in \eqref{adn1}, then we obtain a standard 
terminal value problem for an explicit linear system of ordinary differential equations in
the unknowns $p_n^1,\ldots,p_n^n$ in which all of the coefficient functions belong to
$L^2(0,T)$. By Carath\'eodory's theorem (applied backward in time), it has a unique solution in  
$H^1(0,T;\erren)$ that specifies $p_n$. At the same time, inserting $e_i$, $1\le i\le n$, in
\eqref{adn2}, we recover $(q_n^1,\ldots,q_n^n)\in L^2(0,T;\erren)$, which in turn specifies
$q_n$. Hence, the system \eqref{adn1}--\eqref{adn3} has a unique solution $(p_n,q_n)
\in H^1(0,T;V_n)\times L^2(0,T;V_n)$. 

We now derive a priori estimates. First, we insert $v=p_n$ in \eqref{adn1} and $v=q_n$ in \eqref{adn2},
add the results, and integrate over $(t,T]$ where $t\in [0,T)$. \pier{Noting a cancellation of two terms and recalling} the notation introduced
in \eqref{defQt}, we then arrive at the identity
\begin{align}
&\frac 12\,\|p_n(t)\|^2 \,+\intet |p_n|^2 \,+\intet |q_n|^2\nonumber\\
&=\,\frac 12\,\|p_n(T)\|^2\,-\intet f''(\phistar)\,p_n\,q_n\,+\intet g_1\,p_n\,-\intet g_3\,q_n\,.
\label{adn4}
\end{align} 
Applying Young's inequality, we see that the three integrals on the \rhs\ are bounded by an expression
of the form
$$  
c\,+\,\frac 12 \intet|q_n|^2\,+\,c\intet|p_n|^2, 
$$
and the first summand is bounded since
$\|p_n(T)\|^2=(p_n(T),g_2)\le \|p_n(T)\|\,\|g_2\|$, that is, $\|p_n(T)\|\le \|g_2\|$. 
Therefore, applying Gronwall's lemma backward in time, we obtain from \eqref{adn4} the estimate
\begin{align}
&\|p_n\|_{L^\infty(0,T;H)}\,+\,\|q_n\|_{L^2(0,T;H)}\,\le\,c\quad\forall\,n\in\enne .
\label{adn5}
\end{align}

In the second estimate, we take $v=-\Delta p_n$ in \eqref{adn1} and $v=-\Delta q_n$ in \eqref{adn2}
to obtain the identity
\begin{align}
&\frac 12 \|\nabla p_n(t)\|^2\,+\intet|\nabla p_n(t)|^2\,+\intet|\nabla q_n|^2\non\\
&=\,\frac 12 \|\nabla p_n(T)\|^2+\,\intet\nabla g_1\cdot\nabla p_n\,-\intet\nabla g_3\cdot\nabla q_n
\,-\intet\nabla(f''(\phistar)\,q_n)\cdot\nabla p_n\non\\
&:= I_1+I_2+I_3+I_4,
\label{Bobby}
\end{align}
with obvious meaning. Since $g_1, g_3\in L^2(0,T;V)$, Young's inequality implies that for every 
$\delta>0$ (which has yet to be chosen) we have that
\Beq
I_2+I_3\,\le\,\delta\intet\bigl(|\nabla p_n|^2+|\nabla q_n|^2\bigr)\,+\,\nju{c_\delta}.
\label{Solo}
\Eeq
Moreover, we have, as seen above, $\|p_n(T)\|\,\le \|g_2\|$. Also, by the same token, and since
$\Delta p_n(T)\in \pier{{}V_n{}}$,
$$\|\nabla p_n(T)\|^2=(p_n(T),-\Delta p_n(T))=(g_2,-\Delta p_n(T))=(\nabla g_2,\nabla p_n(T)),$$
which implies that $\|\nabla p_n(T)\|\le \|g_2\|_V$ and thus $\,I_1\,\le\,c$.
  
Finally, using the Young and H\"older inequalities, as well as \eqref{finale}, we obtain that
\begin{align}
I_4\,&=\,-\intet f''(\phistar)\,\nabla q_n \cdot \nabla p_n -\intet f'''(\phistar)\,q_n\,
\nabla\phistar\cdot\nabla p_n\non\\
&\le\,\delta\intet|\nabla q_n|^2\,+\,c_\delta \intet|\nabla p_n|^2
\,+\,c\int_t^T\|q_n(s)\|_4\,\|\nabla\phistar(s)\|_4\,\|\nabla p_n(s)\|_2\,ds\non\\
&\le\,2\delta\intet|\nabla q_n|^2 
\gianni{{}+\delta\intet |q_n|^2}
+ c_\delta\intet \betti{|\nabla p_n|^2}.
\label{sepiangi}
\end{align}
Combining \eqref{Bobby}--\eqref{sepiangi}, \nju{choosing $\,\delta=1/4$,} and using \eqref{adn5} and Gronwall's lemma,
we have thus shown that
\begin{equation}
\|p_n\|_{L^\infty(0,T;V)} \,+\,\|q_n\|_{L^2(0,T;V)}\,\le \,c\quad\forall\,n\in\enne. 
\label{amore}
\end{equation}
But then we may insert $v=-\Delta p_n(t)$ in \eqref{adn2}, and Young's inequality and elliptic
regularity yield that
\begin{equation}
\|\gianni{p_n}\|_{L^2(0,T;W)} \,\le\,c\quad\forall\,n\in\enne.
\label{iopiango}
\end{equation}
\gianni{In our last step, we prove that}
\Beq
\|\dt\gianni{p_n}\|_{L^2(0,T;V^*)} \,\le\,c\quad\forall\,n\in\enne.
\label{conte}
\Eeq
\gianni{%
To this end, we recall that $\Vn$ contains the constant functions, so that \eqref{adn1} yields
\Beq
  - \frac d{dt} \, \barpn + \barpn 
  = - \overline{f''(\phistar)\qn} + \overline{g_1} \,.
  \label{meanpn}
\Eeq
By taking advantage of \eqref{adn3} and \eqref{adn5} for~$\qn$, we deduce that
\Beq
  \norma\barpn_{H^1(0,T)} \leq c \quad\forall\,n\in\enne.
  \label{stimameanpn}
\Eeq
By combining \eqref{adn1} with \eqref{meanpn}, and setting for brevity
\Beq
  \rho_n := - f''(\phistar)\qn + \overline{f''(\phistar)\qn} + g_1 - \overline{g_1}\,,
  \non
\Eeq
we deduce that
\Bsist
  && - \bigl( \dt(\pn-\barpn) , v \bigr)
  + (\pn-\barpn,v)
  + ( \nabla\qn,\nabla v)
  = (\rho_n,v)
  \non
  \\
  && \quad \hbox{\aet, for every $v\in\Vn$} .
  \non
\Esist
\pier{Recalling \eqref{pier2} and \eqref{defN}, we can test the  
 above identity by $-\calN(\dt(\pn-\barpn))$ and} obtain~that
\Bsist
  && \int_t^T \normaVp{\dt(\pn-\barpn)(s)}^2 \, ds
  + \frac 12 \, \normaVp{(\pn-\barpn)(t)}^2
  \non
  \\
  && = \frac 12 \, \normaVp{(\pn-\barpn)(T)}^2
  \pier{{}- \int_t^T \langle \dt(\pn-\barpn)(s), q_n (s)\rangle ds{}}
  - \intet \rho_n \, \calN(\dt(\pn-\barpn)).
  \non
\Esist
\pier{On account of} \eqref{adn3} and \eqref{amore}, Young's inequality, and the properties of~$\calN$,
we see that the whole \rhs\ is bounded~by
\Beq
  \frac 12 \int_t^T \normaVp{\dt(\pn-\barpn)(s)}^2 \, ds
  + c\,,
  \non
\Eeq
and we we infer that \nju{$\,\{\dt(\pn-\barpn)\}\,$} is bounded in $\L2\Vp$.
Then, \eqref{conte} follows on account of~\eqref{stimameanpn}.}

At this point, we conclude from \eqref{amore}--\eqref{conte} the existence of limit points
$(p,q)$ such that, at least for a subsequence \pier{labelled} again by $n$,
\begin{align}
p_n\to p&\quad\mbox{weakly star in }\,H^1(0,T;V^*)\cap L^\infty(0,T;V)\cap L^2(0,T;W)\,,
\label{seridi}
\\
q_n\to q&\quad\mbox{weakly in }\,L^2(0,T;V)\,.
\label{iorido}
\end{align}
It is then a standard matter to conclude that $(p,q)$ is a solution to \Adjoint. We may allow ourselves
to leave this simple argument and the proof of uniqueness to the reader.
\Edim

\step Proof of Theorem \gianni{\ref{Optimality}}

We insert $v=p(t)$ in  \eqref{primaL} and $v=-q(t)$ in \eqref{secondaL}, add the resulting equations, and integrate over $[0,T]$. Using the properties of the adjoint system \Adjoint, we then obtain that
\begin{align}
0\,&=\,\int_0^T\langle \dt\xi(t),p(t)\rangle\,dt \,+\int_Q \bigl(\xi p+\nabla \eta\cdot\nabla p-hp\bigr)
\,+\int_Q\bigl(-\eta q+f''(\phistar)\xi q +\nabla\xi\cdot\nabla q\bigr)\non\\
&=\,\int_0^T\langle -\dt p(t),\xi(t)\rangle\,dt\,+\int_Q \xi\bigl(p+f''(\phistar)q\bigr)\,+
\int_Q\nabla\xi\cdot\nabla q\,+\int_Q\bigl(-\eta q+\nabla\eta\cdot\nabla p\bigr)\non\\
&\qquad -\int_Q hp\,+\iO\xi(T)p(T)\non\\
&=\,\int_Q\bigl(\alpha_1(\phistar-\phi_Q)\xi\,+\,\alpha_3(\mustar-\mu_Q)\eta\,-\,hp\bigr)
\,+\alpha_2\iO(\phistar(T)-\phi_\Omega)\xi(T)\,.\non
\end{align}
Substitution of this identity with $h=u-\ustar$ in \eqref{vug1} yields the
variational inequality \eqref{optimality}.\QED


\section*{Acknowledgments}
\gianni{This research was supported by the Italian Ministry of Education, 
University and Research~(MIUR): Dipartimenti di Eccellenza Program (2018--2022) 
-- Dept.~of Mathematics ``F.~Casorati'', University of Pavia. 
In addition, {PC and ER gratefully mention} some other support 
from the GNAMPA (Gruppo Nazionale per l'Analisi Matematica, 
la Probabilit\`a e le loro Applicazioni) of INdAM (Isti\-tuto 
Nazionale di Alta Matematica).}

\vspace{3truemm}

\Begin{thebibliography}{10}

\betti{%
\bibitem{BO}
M. Bahiana, Y. Oono, 
Cell dynamical system approach to block copolymers, 
{\it Phys. Rev. A} \textbf{41} (1990), 6763-6771.}

\bibitem{Barbu}
V. Barbu,
``Nonlinear Differential Equations of Monotone Type in Banach Spaces'',
Springer,
London, New York, 2010.

\betti{%
\bibitem{BEG}A.L. Bertozzi, S. Esedoglu, A. Gillette,
Inpainting of binary images using the Cahn--Hilliard equation,
{\it IEEE Trans. Image Process.} \textbf{16} (2007), 285-291.}

\betti{%
\bibitem{BE} J.F. Blowey, C.M. Elliott, 
The Cahn--Hilliard gradient theory for phase separation with nonsmooth free energy. I. Mathematical analysis, 
{\it European J. Appl. Math.} {\bf 2} (1991), 233-280.}

\betti{%
\bibitem{BGM}
S. Bosia, M. Grasselli, A. Miranville, 
On the long-time behavior of a 2D hydrodynamic model for chemically reacting binary fluid mixtures, 
{\it Math. Methods Appl. Sci.} {\bf 37} (2014), 726-743}

\bibitem{Brezis}
H. Brezis,
``Op\'erateurs Maximaux Monotones et Semi-groupes de Contractions
dans les Espaces de Hilbert'',
North-Holland Math. Stud.
{\bf 5},
North-Holland,
Amsterdam,
1973.

\betti{%
\bibitem{C}
J.W. Cahn, On spinodal decomposition, 
{\it Acta Metall.} \textbf{9} (1961), 795-801.}

\betti{%
\bibitem{CH}
J.W. Cahn, J.E. Hilliard, 
Free Energy of a Nonuniform System. I. Interfacial Free Energy, 
{\it J. Chem. Phys.} \textbf{28} (1958), 258-267.}

\betti{%
\bibitem{CRW}
C. Cavaterra, E. Rocca, H. Wu, 
Long-time dynamics and optimal control of a diffuse interface model for tumor growth,
{\it Appl. Math. Optim.} {\bf 83} (2021), 739-787.}

\pier{%
\bibitem{CGH}
P. Colli, G. Gilardi, D. Hilhorst,
On a Cahn--Hilliard type phase field system related to tumor growth. {\em Discret. Cont. Dyn. Syst.} {\bf 35} (2015), 2423-2442.   
\bibitem{CGRS1}
P. Colli, G. Gilardi, E. Rocca, J. Sprekels, 
Vanishing viscosities and error estimate for a Cahn--Hilliard
type phase field system related to tumor growth. {\em Nonlinear Anal. Real World Appl.} {\bf 26} (2015), 93-108.
\bibitem{CGRS2}
P. Colli, G. Gilardi, E. Rocca, J. Sprekels, 
Asymptotic analyses and  error estimates for a Cahn--Hilliard
type phase field system modelling tumor growth. {\em Discret. Contin. Dyn. Syst. Ser. S} {\bf 10} (2017), 37-54.}

\betti{%
\bibitem{CGRS}
P. Colli, G. Gilardi, E. Rocca, J. Sprekels, 
Optimal distributed control of a diffuse interface model of tumor growth,
{\it Nonlinearity}  {\bf 30} (2017), 2518-2546.}

\pier{%
\bibitem{CGS-2014}
P. Colli, G. Gilardi, J. Sprekels,
On the Cahn--Hilliard equation with dynamic boundary conditions
and a dominating boundary potential, {\it J. Math. Anal. Appl.}
{\bf 419} (2014), 972-994.
\bibitem{CGS-2015}
P. Colli, G. Gilardi, J. Sprekels,
A boundary control problem for the pure 
Cahn--Hilliard equation with dynamic boundary conditions,
{\it Adv. Nonlinear Anal.} {\bf 4} (2015), 311-325.
\bibitem{CGS-2018a}
P. Colli, G. Gilardi, J. Sprekels,
On a Cahn--Hilliard system with convection 
and dynamic boundary conditions, 
{\it Ann. Mat. Pura Appl. (4)} {\bf 197} (2018), 1445-1475.
\bibitem{CGS-2018b}
P. Colli, G. Gilardi, J. Sprekels,
Optimal velocity control of a viscous Cahn--Hilliard system with convection and dynamic boundary conditions,
{\em SIAM J. Control Optim.} {\bf 56} (2018), 1665-1691.
\bibitem{CGS-2019a}
P. Colli, G. Gilardi, J. Sprekels,
A distributed control problem for a fractional tumor growth model, 
{\em Mathematics} {\bf 7} (2019), 792.
\bibitem{CGS-2019b}
P. Colli, G. Gilardi, J. Sprekels,
Optimal velocity control of a convective Cahn--Hilliard system with double obstacles and 
dynamic boundary conditions:\ \pier{a} `deep quench' approach,
{\em J. Convex Anal.} {\bf 26} (2019), 485-514. 
\bibitem{CGS-2019c}
P. Colli, G. Gilardi, J. Sprekels,
Recent results on well-posedness and optimal control for a class of 
generalized fractional 
Cahn--Hilliard systems, 
{\em Control \pier{Cybernet.}} {\bf 48} (2019), 153-197. 
\bibitem{CGS-2020}
P. Colli, G. Gilardi, J. Sprekels,
Optimal distributed control of a generalized fractional Cahn--Hilliard system,
{\it Appl. Math. Optim.} {\bf 82} (2020), 551-589.
\bibitem {CGS-2021} 
P. Colli, G. Gilardi, J. Sprekels,
Deep quench approximation and optimal control of general Cahn--Hilliard 
systems with fractional operators and double obstacle potentials, 
{\em Discrete Contin. Dyn. Syst. Ser. S} {\bf 14} (2021), 243-271.%
} 

\pier{%
\bibitem{CS-2021}
P. Colli, A. Signori, 
Boundary control problem and optimality conditions 
for the Cahn--Hilliard equation with dynamic boundary conditions, 
{\em Internat. J. Control} {\bf 94} (2021), 1852-1869.}

\betti{%
\bibitem{FGR}
S.~Frigeri, M.~Grasselli, E.~Rocca,
On a diffuse interface model of tumor growth,
{\it European J. Appl. Math.} {\bf 26} (2015), 215-243.}

\pier{\bibitem{GARL_1}
H. Garcke, K.F. Lam,
Well-posedness of a Cahn--Hilliard system modelling tumour
growth with chemotaxis and active transport,
{\it European. J. Appl. Math.} {\bf 28} (2017), 284--316.
\bibitem{GARL_3}
H. Garcke, K.F. Lam,
Global weak solutions and asymptotic limits 
of a Cahn--Hilliard--Darcy system modelling tumour growth,
{\it AIMS Mathematics} {\bf 1} (2016), 318--360.}

\betti{%
\bibitem{GLR}
H. Garcke, K.-F. Lam, E. Rocca,
Optimal control of treatment time in a diffuse interface model of tumor growth,
{\it Appl. Math. Optim.} {\bf 78} (2018), 495-544.}

\bibitem{GiMiSchi} 
G. Gilardi, A. Miranville, G. Schimperna,
On the Cahn--Hilliard equation with irregular potentials and dynamic boundary conditions,
{\it Commun. Pure Appl. Anal.} {\bf 8} (2009), 881-912. 

\bibitem{GR3}
G. Gilardi, E. Rocca,
Well posedness and long time behaviour
for a singular phase field system of conserved type,
{\it IMA J. Appl. Math.} {\bf 72} (2007), 498-530.

\pier{%
\bibitem{GS-2019}
G. Gilardi, J. Sprekels,
Asymptotic limits and optimal control for the Cahn--Hilliard system with convection 
and dynamic boundary conditions, \emph{Nonlinear Anal.} {\bf 178} (2019), 1-31.}

\bibitem{GGM}
A. Giorgini, M. Grasselli, A. Miranville,
The Cahn--Hilliard--Oono equation with singular potential,
{\it Math. Models Methods Appl. Sci.}
{\bf 27} (2017), 2485-2510.

\betti{%
\bibitem{GLRS}
A. Giorgini, K.-F. Lam, E. Rocca, G. Schimperna,
On the existence of strong solutions to the Cahn--Hilliard--Darcy 
system with mass source, \pier{preprint arXiv:2009.13344 [math.AP] (2020), pp.~1-30.}}


\betti{%
\bibitem{KS}
E. Khain, L.M. Sander,
Generalized Cahn--Hilliard equation for biological applications, 
{\it Phys. Rev. E} \textbf{77} (2008), 051129-1-051129-7.}

\betti{%
\bibitem{MR}
S. Melchionna, E. Rocca,
On a nonlocal Cahn--Hilliard equation with a reaction term, 
{\it Adv. Math. Sci. Appl.} {\bf 24} (2014), 461-497.}

\betti{%
\bibitem{M}
A. Miranville, Asymptotic behavior of the Cahn--Hilliard--Oono equation, 
{\it J. Appl. Anal. Comput.} {\bf 1} (2011), 523-536.}

\betti{%
\bibitem{MRS} 
A. Miranville, E. Rocca, G. Schimperna, 
On the long time behavior of a tumor growth model, 
{\it J. Differential Equations} {\bf 67} (2019),  2616-2642.}

\betti{%
\bibitem{MT}
A. Miranville, R. Temam,
On the Cahn--Hilliard--Oono--Navier--Stokes equations with singular potentials, 
{\it Appl. Anal.} {\bf 95} (2016), 2609-2624.}

\bibitem{MiZe}
A. Miranville, S. Zelik, 
Robust exponential attractors for Cahn--Hilliard type equations with singular potentials, 
{\it Math. Methods Appl. Sci.} {\bf 27} (2004), 545-582.

\pier{\bibitem{Moser}
J. Moser,  
A sharp form of an inequality by N. Trudinger.
{\it Indiana Univ. Math. J.}  {\bf 20}  (1970/71), 1077-1092.}

\betti{%
\bibitem{oonopuri87}
Y. Oono, S. Puri, 
Computationally efficient modeling of ordering of quenched phases, 
{\it Phys. Rev. Lett.} {\bf 58} (1987), 836-839.}

\betti{%
\bibitem{oonopuri88I}
Y. Oono, S. Puri, 
Study of phase-separation dynamics by use of cell dynamical systems. I. Modeling,
{\it  Phys. Rev. A} {\bf 38} (1988), 434-\pier{453}.}

\betti{%
\bibitem{oonopuri88II}
\pier{S. Puri, Y. Oono,} 
Study of phase-separation dynamics by use of cell dynamical systems. II. Two-dimensional demonstrations, 
{\it Phys. Rev. A} {\bf 38} (1988), 1542-\pier{1565}.}

\pier{\bibitem{S_b}
A. Signori,
Optimal treatment for a phase field system of Cahn--Hilliard 
type modeling tumor growth by asymptotic scheme, 
{\it Math. Control Relat. Fields} {\bf 10}  (2020) 305-331.
\bibitem{S_a}
A. Signori,
Vanishing parameter for an optimal control problem modeling tumor growth,
{\it Asymptot. Anal.} {\bf 117} (2020) 43-66.}

\bibitem{Simon}
J. Simon,
{Compact sets in the space $L^p(0,T; B)$},
{\it Ann. Mat. Pura Appl.~(4)\/} 
{\bf 146} (1987), 65-96.

\End{thebibliography}

\End{document}
